\documentclass[12pt]{article}
\usepackage{amsmath}
\usepackage{amsthm}
\usepackage{amsfonts}
\usepackage{amssymb}
\newcommand{\be}{\begin{equation}}
\newcommand{\ee}{\end{equation}}
\newcommand{\bea}{\begin{eqnarray}}
\newcommand{\eea}{\end{eqnarray}}
\newcommand{\bean}{\begin{eqnarray*}}
\newcommand{\eean}{\end{eqnarray*}}
\newcommand{\brray}{\begin{array}}
\newcommand{\erray}{\end{array}}
\newcommand{\ben}{\begin{equation}{nonumber}}
\newcommand{\een}{\end{equation}{nonumber}}

\newtheorem{dfn}{Definition}[section]
\newtheorem{thm}[dfn]{Theorem}
\newtheorem{lema}[dfn]{Lemma}
\newtheorem{pro}[dfn]{Proposition}
\newtheorem{coro}[dfn]{Corollary}
\newtheorem{xmpl}[dfn]{Example}
\newtheorem{rmrk}[dfn]{Remark}

\newcommand{\bdfn}{\begin{dfn}}
\newcommand{\bthm}{\begin{thm}}
\newcommand{\blema}{\begin{lema}}
\newcommand{\bpro}{\begin{pro}}
\newcommand{\bcoro}{\begin{coro}}
\newcommand{\bxmpl}{\begin{xmpl}}
\newcommand{\brmrk}{\begin{rmrk}}

\newcommand{\edfn}{\end{dfn}}
\newcommand{\ethm}{\end{thm}}
\newcommand{\elema}{\end{lema}}
\newcommand{\epro}{\end{pro}}
\newcommand{\ecoro}{\end{coro}}
\newcommand{\exmpl}{\end{xmpl}}
\newcommand{\ermrk}{\end{rmrk}}

\newcommand{\half}{\frac{1}{2}}

\newcommand{\B}{\mathcal B}

\newcommand{\EXP}{\textbf{e}}

\numberwithin{equation}{section}
\begin{document}
\begin{center}
{\bf{\large Unitary processes with independent increments and
representations of
  Hilbert tensor algebras}}\\

\vspace{0.2in} {\large  Lingaraj Sahu {\footnote { Stat-Math Unit,
Indian Statistical Institute, Bangalore Centre , $8^{th}$ Mile,
Mysore Road,
 Bangalore-59, India.
 E-mail~:~lingaraj@isibang.ac.in }},~
Michael Sch\"{u}rmann {\footnote { Institut f\"{u}r Mathematik und
Informatik, F.-L.-Jahn-Strasse  15a, D-17487 Greifswald, Germany.
E-mail: schurman@uni-greifswald.de }}}
\\

and\\
{\large  Kalyan B. Sinha {\footnote { Jawaharlal Nehru Centre for
Advanced Scientific Research, Jakkur, Bangalore-64,India } }
{\footnote { Department of Mathematics, Indian Institute of Science,
Bangalore-12, India.\\ E-mail: kbs\_jaya@yahoo.co.in}}}

\end{center}

\begin{abstract}   The aim of this article is  to  characterize  unitary
increment process  by  a quantum stochastic integral representation
on symmetric Fock space. Under  certain assumptions we have proved
its unitary equivalence to a Hudson-Parthasarathy   flow.
\end{abstract}

\section{Introduction}
In the framework of the theory of quantum stochastic calculus
developed by pioneering work of Hudson and Parthasarathy \cite{hp1},
quantum stochastic differential equations (qsde) of the form
 \be
\label{hpeqn00} dV_t=\sum_{\mu,\nu\ge 0}  V_t L_\nu^\mu
\Lambda_\mu^\nu(dt),~ V_0=1_{\mathbf h \otimes \Gamma}, \ee (where
the coefficients $L_\nu^\mu~: \mu,~\nu~\ge 0$ are operators in the
initial Hilbert space  $\mathbf h$  and $\Lambda_\mu^\nu$ are
fundamental processes in the symmetric Fock space  $\Gamma=
\Gamma_{sym}(L^2(\mathbb R_+, \mathbf k))$ with respect to  a fixed
orthonormal basis  (in short `ONB') $\{E_j: j\ge 1\}$ of  the noise
Hilbert space $\mathbf k$ ) have  been formulated  and conditions
for existence and uniqueness  of a solution  $\{V_t\}$ are  studied
by Hudson and Parthasarathy  and many other authors. In particular
when the coefficients $L_\nu^\mu~: \mu,~\nu~\ge 0$ are bounded
operators satisfying  some conditions   it is  observed that the
solution  $ \{V_t:t\ge 0\}$ is a unitary process.

In \cite{hl}, using integral representation of regular quantum
martingales in symmetric Fock space \cite{ps}, the authors show that
any covariant Fock  adapted unitary evolution (with norm-continuous
expectation semigroup)  $\{V_{s,t} : 0\le s\le t<\infty\}$ satisfies
a quantum stochastic differential equation (\ref{hpeqn00}) with
constant coefficients $L_\nu^\mu\in \B(\mathbf h).$  For situations
where the
 expectation semigroup is not norm continuous, the  characterization
  problem is discussed   in \cite{j,ajl}.
In \cite{lw1,lw2}, by extended semigroup methods,  Lindsay   and
Wills have studied such problems for Fock adapted contractive
operator cocycles and completely positive cocyles.

In this  article  we are  interested in the characterization of
unitary evolutions with stationary and  independent increments   on
$\mathbf h \otimes \mathcal H,$ where $ \mathbf h$  and $  \mathcal
H $  are separable Hilbert spaces. In \cite{s1,s2}, by a
co-algebraic treatment, the author has proved that any weakly
continuous unitary stationary independent increment process on
$\mathbf h\otimes \mathcal H, \mathbf h$ finite dimensional,   is
unitarily equivalent to a Hudson-Parthasarathy flow with constant
operator coefficients; see also \cite{LiSk1,LiSk2}. In this present
paper we treat the case of a unitary stationary independent
increment process on $\mathbf h\otimes \mathcal H, \mathbf h$  not
necessarily finite dimensional, with norm-continuous expectation
semigroup. By a GNS type construction we are able to get the noise
space $\mathbf k$ and the bounded operator coefficients $L_\nu^\mu$
such that the Hudson-Parthasarathy flow equation (\ref{hpeqn00})
admits a unique unitary  solution   and is unitarily equivalent to
the unitary process we started with.

\section {Notation and Preliminaries}
 We assume that  all the Hilbert spaces appearing in this article are complex
 separable with inner product  anti-linear
 in the first variable. For any Hilbert spaces $\mathcal H,\mathcal K ~ ~ \mathcal
 B(\mathcal H,\mathcal K)$  and
 $\mathcal B_1( \mathcal H)$ denote the  Banach space  of bounded linear
operators from  $\mathcal H$ to
  $\mathcal K$ and trace class operators on $\mathcal H$
respectively. For a linear map (not necessarily bounded ) $T$  we
write its domain as  $\mathcal D(T).$ We denote the trace on
$\mathcal B_1( \mathcal H)$ by $Tr_{\mathcal H}$ or simply $Tr.$ The
 von Neumann algebra of bounded linear operators on  $\mathcal H$ is denoted by $B(\mathcal H).$
 The  Banach space   $\mathcal B_1( \mathcal H, \mathcal K)\equiv \{ \rho\in  \mathcal B( \mathcal H, \mathcal K)
 :  |\rho| :=\sqrt{\rho^*\rho }\in \mathcal B_1( \mathcal H) \}$  with norm (Ref. Page no. 47 in   \cite{gelshi})
\[\|\rho\|_1=\|  ~|\rho|~\|_{\mathcal B_1( \mathcal H)}= \sup\{ \sum_{k\ge 1} |\langle \phi_k, \rho \psi_k \rangle |
:\{\phi_k\},\{\psi_k\}~are ~ONB~of~ \mathcal K ~and ~\mathcal H
~resp. \}\] is the predual of
 $\mathcal B( \mathcal K, \mathcal H).$  For an element $ x\in  \mathcal B( \mathcal K, \mathcal H),$ ~
 $  \mathcal B_1( \mathcal H, \mathcal K) \ni \rho \mapsto Tr_{\mathcal H} (x\rho)$ defines an
  element  of the  dual
  Banach space  $\mathcal B_1( \mathcal H, \mathcal K)^*.$
 For a linear map $T$ on the Banach space $\mathcal B_1( \mathcal H, \mathcal K)$ the adjoint $T^*$ on the dual
  $\mathcal B( \mathcal K, \mathcal H)$ is given by
 $ Tr_{\mathcal H} (T^* (x)\rho):=Tr_{\mathcal H} (x  T(\rho)),~ \forall  x\in  \mathcal B( \mathcal K, \mathcal H),
 ~\rho \in \mathcal B_1( \mathcal H, \mathcal K). $\\

\noindent
 For any   $\xi \in \mathcal H \otimes
\mathcal K, h \in \mathcal H$  the map \[ \mathcal K \ni k
\mapsto\langle \xi, h\otimes k
 \rangle\] defines  a bounded  linear  functional on $\mathcal K$
  and thus  by  Riesz's  theorem
there exists  a unique vector  $\langle \langle \xi, h\rangle
\rangle$ in $ \mathcal K$  such that
 \be \label{partinn}
\langle~\langle \langle \xi, h\rangle \rangle,~k\rangle= \langle
\xi, h\otimes k\rangle, \forall k\in \mathcal K.\ee In other words
$\langle \langle \xi, h\rangle \rangle=F_h^*\xi$ where  $F_h\in
\B(\mathcal K,\mathcal H \otimes \mathcal K)$ is given by $F_h k=h
\otimes k.$

 Let $\mathbf h$ and $\mathcal H$ be two Hilbert spaces with some
  orthonormal bases $\{e_j:j \ge 1\}$ and
 $\{\zeta_j:j\ge 1\}$ respectively.
For $A\in \mathcal B( \mathbf h \otimes \mathcal H)$ and $u,v\in
\mathbf h$ we define a  linear operator $A(u,v)\in \mathcal B(
\mathcal H)$ by
 \[\langle \xi_1,A(u,v) \xi_2\rangle=\langle u\otimes \xi_1 ,A ~v\otimes \xi_2 \rangle,~\forall \xi_1,\xi_2 \in \mathcal
H\]  and read off the following properties:
 \blema \label{Auv}
 Let $A,B \in \mathcal B( \mathbf h \otimes \mathcal H)$ then  for any $u,v,u_i$ and $v_i, i=1,2$ in $\mathbf h$
 \begin{description}
 \item  [(i)] $A(u,v)\in \mathcal B(\mathcal H)$ with $\|A(u,v)\| \le  \|A\|~\|u\|~ \|v\|$ and  $A(u,v)^*=A^*(v,u).$
 \item  [(ii)] $\mathbf h \times \mathbf h \mapsto A(\cdot~,\cdot)$ is
 $1-1,$ i.e. if  $A(u,v)=B(u,v),~\forall u,v \in \mathbf h $ then $A=B.$
 \item  [(iii)]  $A(u_1,v_1)B(u_2,v_2)=[A(|v_1><u_2|\otimes 1_{\mathcal H})B](u_1,v_2)$
 \item  [(iv)] $AB(u,v)=\sum_{j \ge 1} A(u,e_j) B(e_j,v)$ (strongly)
 \item  [(v)]  $0\le A(u,v)^* A(u,v) \le \|u\|^2 A^*A(v,v)$
 \item  [(vi)] $ \langle A(u,v)\xi_1, B(p,w)\xi_2 \rangle= \sum_{j\ge 1}\langle
 p \otimes \zeta_j, [B( |w><v|\otimes |\xi_2><\xi_1|)A^* u\otimes \zeta_j
\rangle\\
~~~~~~~~~~~~~~~~~~~~~~~~~~=\langle v\otimes \xi_1,~ [A^*(
|u><p|\otimes 1_{\mathcal H})B w \otimes \xi_2 \rangle .$
 \end{description}
 \elema
 \begin{proof}
 We are omitting the proof of (i),(ii).\\
\noindent {\bf (iii)} For any $\xi, \zeta\in \mathcal H$ we have
\bean
 && \langle \xi,A(u_1,v_1)B(u_2,v_2)\zeta \rangle=\langle u_1\otimes \xi,A v_1
  \otimes B(u_2,v_2)\zeta \rangle=\langle A^* u_1\otimes \xi, v_1 \otimes B(u_2,v_2)
  \zeta \rangle\\
&&=\sum_{n\ge 1}\langle A^* u_1\otimes \xi, v_1 \otimes
\zeta_n\rangle
\langle \zeta_n, B(u_2,v_2)\zeta \rangle\\
&&=\sum_{n\ge 1}\langle A^* u_1\otimes \xi, v_1 \otimes
\zeta_n\rangle
\langle u_2 \otimes \zeta_n, B v_2\otimes \zeta \rangle\\
&&=\sum_{n\ge 1}\langle A^* u_1\otimes \xi, (|v_1><u_2| \otimes |\zeta_n><\zeta_n|)
 B v_2\otimes \zeta \rangle\\
&&=\langle  u_1\otimes \xi, A (|v_1><u_2| \otimes 1_{\mathcal H} ) B
v_2\otimes \zeta \rangle. \eean Thus it follows that
\[A(u_1,v_1)B(u_2,v_2)=[A(|v_1><u_2|\otimes 1_{\mathcal H})B](u_1,v_2).\]
\noindent {\bf (iv)} By  part(iii) \bean
&&\sum_{j= 1}^N \|A(e_j,u) \xi\|^2\\
&&=\sum_{j= 1}^N \langle \xi, A^*(u,e_j) A(e_j,u) \xi \rangle\\
&&= \langle \xi, [A^*(P_N\otimes 1_{\mathcal H})A](u,u) \xi \rangle,
\eean where $P_N$ is the finite rank projection $\sum_{j= 1}^N
|e_j><e_j|$ on $\mathbf h.$ Since $\{[A^*(P_N\otimes 1_{\mathcal
H})A](u,u)\}$ is an increasing  sequence of positive operators  and
$0\le P_N\otimes 1_{\mathcal H}$  converges strongly to  $1_{\mathbf
h \otimes \mathcal H}$ as $N$ tends to $\infty,$ $[A^*(P_N\otimes
1_{\mathcal H})A](u,u)$  converges strongly to  $[A^*A](u,u)$  as
$N$ tends to $\infty.$ Thus
\[\lim_{N\rightarrow \infty}\sum_{j= 1}^N \|A(e_j,u) \xi\|^2=\langle \xi,
 [A^*A](u,u) \xi \rangle\] and
\[\sum_{j= 1}^N \|A(e_j,u) \xi\|^2\le \| A ~u \otimes  \xi\|^2\le
\| A\|^2  \|u\|^2 \| \xi\|^2 , \forall N\ge 1.\] Now let us consider
the following, for $\xi, \zeta \in  \mathcal H$ \bean
&& |\langle \xi, \sum_{j = 1}^N A(u,e_j) B(e_j,v) \zeta \rangle|^2 =|\sum_{j = 1}^N
\langle A^*(e_j,u) \xi,   B(e_j,v) \zeta \rangle|^2\\
&&\le \sum_{j = 1}^N \|A^*(e_j,u) \xi \|^2 \sum_{j = 1}^N \| B(e_j,v) \zeta \|^2\\
&&\le \| A\|^2  \|u\|^2 \| \xi\|^2 \| B\|^2  \|v\|^2 \| \zeta\|^2.
\eean So
\[ |\langle \xi, \sum_{j = 1}^N A(u,e_j) B(e_j,v) \zeta \rangle| \le \| A\| \| B\|
 \|u\| \|v\| \|\xi\| \| \zeta\| \] and strong convergence of
$\sum_{j \ge 1} A(u,e_j) B(e_j,v)$ follows.\\

 \noindent (v)
 We have
 \bean
 && \langle \xi,A(u,v)^* A(u,v)\xi \rangle
 =\sum_{j\ge 1}\langle \xi, A^* (v,u) \zeta_j\rangle \langle \zeta_j, A (u,v)
 \xi\rangle\\
 && =\sum_{j\ge 1}\langle v\otimes \xi, A^* u \otimes  \zeta_j\rangle \langle
 u \otimes \zeta_j,
  A v\otimes \xi\rangle\\
  && =\langle v\otimes \xi, A^* \{| u ><u| \otimes \sum_{j\ge 1}|
  \zeta_j><\zeta_j|\}
  A v\otimes \xi\rangle.
  \eean
  Since $\sum_{j\ge 1}|
  \zeta_j><\zeta_j|$ converges strongly to the  identity operator
  \[\langle \xi,A(u,v)^* A(u,v)\xi \rangle \le \|u\|^2 \langle v\otimes \xi, A^*
  A v\otimes \xi\rangle\] and this proves the result.

\noindent {\bf (vi)}We have \bean
&&\langle A(u,v)\xi_1, B(p,w)\xi_2 \rangle\\
&&=\sum_{j\ge 1}\langle A(u,v)\xi_1 , \zeta_j \rangle \langle   \zeta_j,  B(p,w)\xi_2
 \rangle\\
&&=\sum_{j\ge 1} \langle   A v \otimes  \xi_1, u \otimes  \zeta_j \rangle \langle  p
 \otimes  \zeta_j,  B w \otimes \xi_2 \rangle\\
&&=\sum_{j\ge 1} \langle B^* p \otimes  \zeta_j,   w \otimes \xi_2
\rangle
\langle   v \otimes  \xi_1, A^* u \otimes  \zeta_j \rangle \\
&&=\sum_{j\ge 1} \langle p \otimes  \zeta_j,  B(| w><v| \otimes
|\xi_2 >< \xi_1|) A^* u \otimes  \zeta_j \rangle. \eean This proves
the  first  part of (vi), the other part follows from \bean
&&\sum_{j\ge 1} \langle  p \otimes  \zeta_j,  B(| w><v| \otimes
|\xi_2 >< \xi_1|)
 A^* u \otimes  \zeta_j \rangle\\
&&= Tr_{\mathbf h \otimes \mathcal H} [( |u><p|\otimes 1_{\mathcal H} )
B(| w><v| \otimes |\xi_2 >< \xi_1|) A^*  ] \\
&&= Tr_{\mathbf h \otimes \mathcal H} [ (| w><v| \otimes |\xi_2 >< \xi_1|)
A^*  ( |u><p|\otimes 1_{\mathcal H} ) B  ] \\
&&=\langle v\otimes \xi_1,~ [A^*( |u><p|\otimes 1_{\mathcal H})B w \otimes \xi_2 \rangle\\
\eean
 \end{proof}

\subsection{Symmetric Fock Space  and Quantum Stochastic Calculus}

Let us briefly recall the
 fundamental integrator processes  of
quantum stochastic calculus  and the flow equation, introduced by
Hudson and Parthasarathy \cite{hp1}. For a Hilbert space $\mathbf k$
let us consider the symmetric Fock space $\Gamma=\Gamma(L^2(\mathbb
R_+,\mathbf k)).$  The exponential vector in  the Fock space,
associated with a vector $f\in L^2(\mathbb R_+, \mathbf k)$ is given
by
 \[\EXP(f)=\bigoplus_{n\ge 0}\frac{1}{\sqrt{n!}}f^{(n)},\]
 where $ f^{(n)}=\underbrace{f \otimes f  \otimes
 \cdots \otimes  f}_{n-copies} $ for $ n>0 $
and by convention  $f^{(0)}=1 .$ The exponential vector
$\textbf{e}(0)  $ is called    the vacuum vector.

\noindent  Let us consider the Hudson-Parthasarathy  (HP) flow
equation on $\mathbf h \otimes \Gamma(L^2(\mathbb R_+, \mathbf k))$:

\be \label{hpeqn} V_{s,t}=1_{\mathbf h \otimes
\Gamma}+\sum_{\mu,\nu\ge 0}  \int_s^t V_{s,\tau} L_\nu^\mu
\Lambda_\mu^\nu(d\tau).\ee Here the coefficients $L_\nu^\mu~:
\mu,~\nu~\ge 0$ are operators in $\mathbf h$  and $\Lambda_\mu^\nu$
are fundamental processes with respect to a fixed orthonormal basis
$\{E_j: j\ge 1\}$ of $\mathbf k:$

 \begin{equation}
 \Lambda_\nu^\mu (t)=\left\{
 \begin{array} {lll} & t ~1_{\mathbf h \otimes
\Gamma} & \mbox{for}
(\mu,\nu)=(0,0) \\
& a(1_{[0,t]}\otimes E_j)&  \mbox{for} (\mu,\nu)=(j,0)
\\
& a^\dag(1_{[0,t]}\otimes E_k)& \mbox{for}(\mu,\nu)=(0,k)
\\
& \Lambda (1_{[0,t]}\otimes |E_k><E_j|)& \mbox{for} (\mu,\nu)=(j,k).
 \end{array}
 \right.
  \end{equation}

\bthm \label{hpflow} \cite{ms1,krp,gs}  Let $H\in \mathcal
B({\mathbf{h}})$ be self-adjoint, $\{L_k,~W_k^j:j,k \ge 1\}$
 be a family of bounded linear operators in ${\mathbf{h}}$ such that
$W=\sum_{j,k\ge 1}W_k^j\otimes |E_j><E_k|$ is an isometry
(respectively  co-isometry) operator in ${\mathbf{h}}\otimes
{\mathbf{k}}$ and for some
 constant $c\ge 0,$
\[\sum_{k\ge 1}\|L_k u\|^2\le c \|u\|^2,\ \forall u\in {\mathbf{h}}.\]
Let  the coefficients $L_\nu^\mu$ be as follows,
\begin{equation}
L_\nu^\mu=\left\{ \begin{array} {lll}
 & iH-\half\sum_{k\ge 1}L_k^*L_k & \mbox{for}\
(\mu,\nu)=(0,0)\\
&L_j &  \mbox{for}\ (\mu,\nu)=(j,0)
\\
& -\sum_{j\ge 1}L_j^*W_k^j &
\mbox{for}\ (\mu,\nu)=(0,k)\\
& W_k^j- \delta_k^j & \mbox{for}\ (\mu,\nu)=(j,k).
 \end{array}
 \right.
  \end{equation}
   Then  there exists a  unique isometry ( respectively
co-isometry) operator valued process $V_{s,t}$ satisfying  \
(\ref{hpeqn}) . \ethm

\section{Hilbert tensor algebra}

For a product vector $\underbar{u}=u_1\otimes u_2\otimes  \cdots
\otimes u_n\in \mathbf h^{\otimes n}$ we shall denote the product
vector $u_n\otimes u_{n-1}\otimes  \cdots \otimes u_1$ by
$\underleftarrow{\mbox{u}}.$ For the null vector in $\mathbf
h^{\otimes n}$  we  shall write $\underbar{0}$. If
$\{f_j\}_{j=1}^\infty$ is an ONB for
 $\mathbf h,$ then we have a product ONB
 $\{f_{\underbar{j}}=f_{j_1}\otimes \cdots \otimes f_{j_n}: \underbar{j}
 =(j_1,j_2,
\cdots, j_n) , j_k \ge 1\}$ for the Hilbert space
$\mathbf h^{\otimes n}.$\\

Consider $\mathbb Z_2=\{0,1\},$ the finite field
 with addition modulo $2.$
 For $n\ge 1,$ let $ \mathbb Z_2^n$ denotes the $n$-fold direct sum
 of
$\mathbb Z_2$  and we write
 $\underbar{0}=(0,0, \cdots, 0)$ and $\underbar{1}=(1,1, \cdots, 1).$  For
$\underline{\epsilon}=(\epsilon_1,\epsilon_2,\cdots, \epsilon_n),
\underline{\epsilon}^\prime=(\epsilon_1^\prime,\epsilon_2^\prime,\cdots,
\epsilon_m^\prime)$ we put
\[\underline{\epsilon}\oplus \underline{\epsilon}^\prime=
(\epsilon_1,\cdots, \epsilon_n,\epsilon_1^\prime,\cdots,
\epsilon_m^\prime)\in \mathbb Z_2^{n+m}  ~\mbox{and we define }~
\underline{\epsilon}^*=\underline{1}+
(\epsilon_n,\epsilon_{n-1},\cdots, \epsilon_1)\in \mathbb Z_2^n.\]

Let $A \in \mathcal B( \mathbf h \otimes \mathcal H),\epsilon \in
\mathbb Z_2=\{0,1\}.$ We define operators  $A^{(\epsilon)}\in
\mathcal B(\mathbf h\otimes \mathcal H)$ by $A^{(\epsilon)}:= A$ if
$\epsilon=0$  and $A^{(\epsilon)}:= A^*$ if $\epsilon=1.$  For $1\le
k\le n,$ we define a unitary exchange map $P_{k,n}:\mathbf
h^{\otimes n} \otimes \mathcal H\rightarrow \mathbf h^{\otimes n}
\otimes \mathcal H  $ by putting
 \[P_{k,n}(  \underbar{u}  \otimes  \xi ) :=u_1\otimes \cdots
  \otimes u_{k-1} \otimes u_{k+1}\cdots  \otimes  u_n \otimes ( u_k \otimes \xi)   \]
    on product vectors.
Let $\underline{\epsilon}=(\epsilon_1,\epsilon_2,\cdots,
\epsilon_n)\in \mathbb Z_2^n.$  Consider the  ampliation of the
operator $A^{(\epsilon_k)}$ in $\mathcal B(\mathbf h^{\otimes n}
\otimes \mathcal H)$ given by
  \[ A^{(n,\epsilon_k)}:=P_{k,n}^* (1_{\mathbf h^{\otimes n-1}} \otimes A^{(\epsilon_k)})P_{k,n}.\]
   Now we define the operator   $A^{(\underline{\epsilon})}:=\prod
_{k=1}^n~A^{(n,\epsilon_k)}:=A^{(1,\epsilon_1)}\cdots
A^{(n,\epsilon_n)}$ in
  $\mathcal B(\mathbf h^{\otimes n} \otimes \mathcal H).$ Please note
  that as
here, through out  this article, the product symbol  $\prod
_{k=1}^n$  stands for  product with order  $1$ to $n.$ For $m\le
  n,$
  we shall write $\underline{\epsilon}^{(m)}=(\epsilon_1,\epsilon_2,\cdots,
  \epsilon_m)$  and consider the operator
$A^{(\underline{\epsilon}^{(m)})}=\prod_{i=1}^m A^{(m,\epsilon_i)}$
in $\mathcal B(\mathbf h^{\otimes m} \otimes \mathcal H)$ \noindent
We have the  following  preliminary observation.
\blema \label{Akn}

 \begin{description}
\item[(i)]
 For product vectors  $\underbar{u}, \underbar{v}\in
\mathbf h^{\otimes n}$
 \[\prod_{i=1}^m A^{(n,\epsilon_i)} (\underbar{u},\underbar{v})
 =\prod_{i=1}^m A^{\epsilon_i}(u_i,v_i)   \prod_{i=m+1}^n \langle u_i,v_i\rangle
      \in  \mathcal B(\mathcal H).\]
\item[(ii)] For $\xi, \zeta\in \mathcal H$
\[ \prod_{i=1}^m A^{(\epsilon_i)}(\xi, \zeta)=A^{(\underline{\epsilon}^{(m)})}
(\xi, \zeta) \otimes 1_{\mathbf h^{\otimes n-m}} \in B(\mathbf
h^{\otimes n}).\]
\item[(iii)] If $A$ is  an isometry (respectively  unitary) then $A^{(n,\epsilon_k)}$ and
$A^{(\underline{\epsilon})}$ are isometries (respectively
unitaries).
 \end{description}
\elema
\noindent The proof is obvious and is omitted.\\
We note that part (i) of this Lemma  in particular gives \be
\label{A-prod}
A^{(\underline{\epsilon})}(\underbar{u},\underbar{v})=\prod_{i=1}^n
A^{(\epsilon_i)}(u_i,v_i)
 \ee

\noindent Let
$M_0:=\{(\underbar{u},\underbar{v},\underline{\epsilon}):~\underbar{u}=\otimes_{i=1}^n
u_i,~\underbar{v}=\otimes_{i=1}^n v_i\in   \mathbf h^{\otimes n},
\underline{\epsilon}=(\epsilon_1,\epsilon_2,\cdots, \epsilon_n)\in
\mathbb Z_2^n, ~ n\ge 1\}.$ In $M_0,$ we introduce an equivalence
relation $`\sim $'~
$:~(\underbar{u},\underbar{v},\underline{\epsilon})\sim
(\underbar{p},\underbar{w},\underline{\epsilon}^\prime)$ if
$\underline{\epsilon}=\underline{\epsilon}^\prime$ and
 $|\underbar{u}><\underbar{v}|=|\underbar{p}>< \underbar{w}|\in
 \B(\mathbf h^{\otimes n}).$
 Expanding the vectors in
term of the ONB $\{e_{\underbar{j}}=e_{j_1}\otimes \cdots \otimes
e_{j_n}: \underbar{j}  =(j_1,j_2, \cdots, j_n) , j_k\ge 1\}$, from
$|\underbar{u}><\underbar{v}|=|\underbar{p}>< \underbar{w}|$ we get
${\underbar{u}}_{\underbar {j}}{\overline{\underbar{v}}}_{\underbar
{k}}
    ={\underbar{p}}_{\underbar {j}}
{\overline{\underbar{w}}}_{\underbar {k}}
    $ for each
multi-indices $\underbar {j},\underbar {k}.$ Thus in particular when
$(u,v,0)\sim (p,w,0), $  for any $\xi_1,\xi_2\in \mathcal H$ we have
\bean && \langle \xi_1, A(u,v) \xi_2 \rangle\\
&& =\sum_{j,k\ge 1} \overline{u_j} v_k \langle  e_j \otimes \xi_1, A
e_k \otimes  \xi_2 \rangle \\
&&=\sum_{j,k\ge 1} \overline{p_j} w_k \langle  e_j \otimes \xi_1, A
e_k \otimes  \xi_2 \rangle\\
&&= \langle \xi_1, A(p,w) \xi_2 \rangle. \eean In fact
$A(u,v)=A(p,w)$ iff $(u,v,0)\sim (p,w,0)$ and more generally
$A^{(\underline{\epsilon})}(\underbar{u},\underbar{v})
=A^{(\underline{\epsilon}^\prime)}(\underbar{p},\underbar{w})$ iff
$(\underbar{u},\underbar{v},\underline{\epsilon})\sim
(\underbar{p},\underbar{w},\underline{\epsilon}^\prime).$ It is easy
to see that $(\underbar{0},\underbar{v},\underline{\epsilon})\sim
(\underbar{u},\underbar{0},\underline{\epsilon})\sim
(\underbar{0},\underbar{0},\underline{\epsilon})$ and we call this
class the $0$ of the quotient set $M_0.$

Let us define  multiplication and involution on $M_0/\sim$ by setting\\
Vector multiplication: ~ $(\underbar{u},\underbar{v},\underline{\epsilon}).
(\underbar{p},\underbar{w},\underline{\epsilon}^\prime)= (\underbar{u}
\otimes \underbar{p},\underbar{v} \otimes \underbar{w},\underline{\epsilon}\oplus
\underline{\epsilon}^\prime) $  and \\
Involution:  $(\underbar{u},\underbar{v},\underline{\epsilon})^*=
 (\underleftarrow{\mbox{v}},\underleftarrow{\mbox{u}},\underline{\epsilon}^*).$\\

 \noindent
 Since $\underleftarrow{(\underbar{u} \otimes \underbar{p} )}=
 p_m\otimes \cdots   \otimes  p_1 \otimes u_n \otimes \dots   \otimes
 u_1= (\underleftarrow{\mbox{p}}\otimes \underleftarrow{\mbox{u}})
 $ and $(\underline{\epsilon}\oplus\underline{\epsilon}^\prime)^*=
  (\underline{\epsilon}^\prime)^* \oplus \underline{\epsilon}^*$
 \bean
&&[(\underbar{u},\underbar{v},\underline{\epsilon}).
(\underbar{p},\underbar{w},\underline{\epsilon}^\prime)]^*=
(\underbar{u}\otimes \underbar{p} ,\underbar{v}\otimes
\underbar{w},\underline{\epsilon}\oplus\underline{\epsilon}^\prime)^*\\
&&=(\underleftarrow{\underbar{v}\otimes \underbar{w}},
\underleftarrow{\underbar{u}\otimes \underbar{p}}
,(\underline{\epsilon}\oplus\underline{\epsilon}^\prime)^*)\\
&&=(\underleftarrow{\mbox{w}}\otimes \underleftarrow{\mbox{v}},
\underleftarrow{\mbox{p}}\otimes \underleftarrow{\mbox{u}} ,
 (\underline{\epsilon}^\prime)^* \oplus \underline{\epsilon}^*)\\
&&=(\underbar{p},\underbar{w},\underline{\epsilon}^\prime)^* .
(\underbar{u},\underbar{v},\underline{\epsilon})^*.\eean It is clear
that $\underline{\epsilon}=\underline{\epsilon}^\prime \implies
\underline{\epsilon}^*=(\underline{\epsilon}^\prime)^*$ and
$|\underbar{u}><\underbar{v}|=|\underbar{p}>< \underbar{w}|$ implies
$|\underleftarrow{\mbox{v}}><\underleftarrow{\mbox{u}}|=|\underleftarrow{\mbox{w}}><
\underleftarrow{\mbox{p}}|.$  Thus $
(\underbar{u},\underbar{v},\underline{\epsilon})\sim
(\underbar{p},\underbar{w},\underline{\epsilon}^\prime)$ implies
$(\underbar{u},\underbar{v},\underline{\epsilon})^*\sim
(\underbar{p},\underbar{w},\underline{\epsilon}^\prime)^*.$
Moreover, $(\underbar{u},\underbar{v},\underline{\epsilon})\sim
(\underbar{u}',\underbar{v}',\underline{\epsilon}^\prime)$  and
$(\underbar{p},\underbar{w},\underline{\alpha})\sim
(\underbar{p}',\underbar{w}',\underline{\alpha}^\prime)$ implies
 $\underline{\epsilon}\oplus
 \underline{\alpha}=\underline{\epsilon}^\prime\oplus
 \underline{\alpha}'$
and $|\underbar{u}\otimes \underbar{p}
><\underbar{v} \otimes \underbar{w}|=|\underbar{u}>< \underbar{v}|
 \otimes |\underbar{p}>< \underbar{w}|=
 |\underbar{u}'>< \underbar{v}'|
 \otimes |\underbar{p}'>< \underbar{w}'|
 =|\underbar{u}'\otimes \underbar{p}'
><\underbar{v}' \otimes \underbar{w}'|.$
 So that the involution and multiplication respect $\sim.$

 Let $M$ be the complex  vector space spanned by  $
M_0/\sim.$ The elements of $M$ are  formal finite linear
combinations of elements of $ M_0/\sim.$   With the above
multiplication and involution $M$ is a $*$-algebra.

\section {Unitary processes with stationary and independent increment }
Let $\{U_{s,t}: 0\le s\le t<\infty\}$ be a family of unitary
operators in $ \mathcal B( \mathbf h \otimes \mathcal H)$ and
$\Omega $ be a fixed unit vector in $\mathcal H.$ We shall also set
$U_t:=U_{0,t}$ for simplicity. As we discussed in the previous
section,  let us consider the family of  operators  $\{
U_{s,t}^{(\epsilon)}\}$ in $ \mathcal B( \mathbf h\otimes \mathcal
H)$  for $\epsilon \in \mathbb Z_2$ given by $ U_{s,t}^{(\epsilon)}=
U_{s,t}$ if $\epsilon=0, U_{s,t}^{(\epsilon)}= U_{s,t}^*$ if
$\epsilon=1.$ Furthermore  for $n\ge 1 , \underline{\epsilon}\in
\mathbb Z_2^n$ fixed, $1\le k\le n, $ we consider the families  of
operators $\{U_{s,t}^{(\epsilon_k)}\}$ and
$\{U_{s,t}^{(\underline{\epsilon})}\}$ in $ \mathcal B( \mathbf
h^{\otimes n} \otimes \mathcal H).$ By Lemma \ref{Akn} we observe
that
\[
U_{s,t}^{(\underline{\epsilon})}(\underbar{u},\underbar{v})
=\prod_{i=1}^n U_{s,t}^{(\epsilon_i)}(u_i,v_i).\]

\noindent For $\underline{\epsilon}=\underline{0} \in \mathbb Z_2^n$
and $1\le k\le n,$ we shall write $U_{s,t}^{(n,k)} $ for the unitary
operator $U_{s,t}^{(n,\epsilon_k)}$ and  $U_{s,t}^{(n)}$ for the
unitary $U_{s,t}^{(\underline{0})}$ on $\mathbf h^{\otimes n}
\otimes \mathcal H.$ For $n\ge 1, \underbar{s}=(s_1,s_2, \cdots,
s_n), \underbar{t}=(t_1,t_2, \cdots, t_n)$ $:~ 0 \le s_1\le t_1\le
s_2\le  \ldots \le  s_n\le  t_n< \infty,$
$\underline{\epsilon}_k=(\alpha_1^{(k)},  \alpha_2^{(k)}, \cdots,
\alpha_{m_k}^{(k)} ) \in \mathbb Z_2^{m_k} : 1\le k\le
n,m=m_1+m_2+\cdots +m_n$
$\underline{\epsilon}=\underline{\epsilon}_1\oplus
\underline{\epsilon}_2 \oplus \cdots \oplus
\underline{\epsilon}_n\in \mathbb Z_2^m, $ we define $
U_{\underbar{s},\underbar{t}}^{(\underline{\epsilon})}\in \mathcal
B(\mathbf h^{\otimes m}\otimes \mathcal H)$ by setting \be
\label{U-underbar-st}
U_{\underbar{s},\underbar{t}}^{(\underline{\epsilon})}:=\prod_{k=1}^n
U_{s_k,t_k}^{(\underline{\epsilon}_k)}.\ee
 Here
$U_{s_k,t_k}^{(\underline{\epsilon}_k)}$ is  looked  upon as an
operator in $\mathcal B(\mathbf h^{\otimes m}\otimes \mathcal H)$ by
ampliation
 and appropriate  tensor flip. So for $\underbar{u}=\otimes_{k=1}^n \underbar{u}_k,
\underbar{v}=\otimes_{k=1}^n \underbar{v}_k\in \mathbf h^{\otimes m}
$ we have
\[U_{\underbar{s},\underbar{t}}^{(\underline{\epsilon})}(\underbar{u},
\underbar{v}) =\prod_{k=1}^n U_{s_k,t_k}^{(\underline{\epsilon}_k)}
(\underbar{u}_k,\underbar{v}_k) .\]  When there can be  no
confusion, for $\underline{\epsilon}=\underline{0}$ we   write
$U_{\underbar{s},\underbar{t}}$ for
$U_{\underbar{s},\underbar{t}}^{(\underline{\epsilon})}.$ For
$a,b\ge 0, \underbar{s}=(s_1,s_2, \cdots, s_n),
\underbar{t}=(t_1,t_2, \cdots, t_n)$ we write $a\le
\underbar{s},\underbar{t}\le b$  if $ a\le s_1\le t_1\le s_2 \le
\ldots \le s_n\le  t_n\le b.$

 \noindent Let us assume the following  properties on the unitary family
$U_{s,t}$  for further discussion to prove  unitary equivalence  of
$U_{s,t}$  with an HP flow.\\ \\
 {\bf Assumption A}
\begin{description}
\item [A1] {\bf (Evolution)} For any $   0\le r\le s\le t<\infty, ~ U_{r,s}U_{s,t}=U_{r,t}.$
\item [A2] {\bf (Independence of increments)} For any ~ $0\le s_i\le t_i<\infty ~:
 ~ i=1,2$ such that $[s_1,t_1)\cap[s_2,t_2)=\varnothing$\\
 (a)~$
U_{s_1,t_1}(u_1,v_1) $ commutes with  $U_{s_2,t_2}(u_2,v_2) $
 and  $U_{s_2,t_2}^*(u_2,v_2) $  for every  $u_i,v_i \in \mathbf
 h.$\\
(b)~  For  $s_1\le \underbar{a},\underbar{b}\le  t_1,~~s_2\le
\underbar{q},\underbar{r}\le  t_2$ and $\underbar{u},\underbar{v}\in
\mathbf h^{\otimes n},~
 \underbar{p},\underbar{w}\in
\mathbf h^{\otimes m},\underline{\epsilon}\in \mathbb
Z_2^n,\underline{\epsilon}^\prime \in \mathbb Z_2^m $
\[\langle \Omega,  U_{\underbar{a},\underbar{b}}^{(\underline{\epsilon})}(\underbar{u},
\underbar{v})
U_{\underbar{q},\underbar{r}}^{(\underline{\epsilon}^\prime)}(\underbar{p},
\underbar{w}) \Omega \rangle = \langle \Omega,
U_{\underbar{a},\underbar{b}}^{(\underline{\epsilon})}(\underbar{u},
\underbar{v})\Omega \rangle
 \langle \Omega,U_{\underbar{q},\underbar{r}}^{(\underline{\epsilon}^\prime)}(\underbar{p},
\underbar{w})  \Omega \rangle.\]
\item [A3]  {\bf (Stationarity)} For any $0\le s\le t<\infty$
 and  $ \underbar{u},\underbar{v}\in \mathbf h^{\otimes n} ,
 \underline{\epsilon}\in \mathbb Z_2^n$
 \[\langle \Omega,  U_{s,t}^{(\underline{\epsilon})}(\underbar{u},\underbar{v}
 ) \Omega \rangle
=  \langle \Omega,
U_{t-s}^{(\underline{\epsilon})}(\underbar{u},\underbar{v}) \Omega
\rangle.\]

\item [Assumption B]
{\bf (Uniform continuity)}\\
 $\lim_{t \rightarrow 0 }~ \sup \{
|\langle \Omega,   (U_t-1)(u,v) \Omega \rangle|: \|u\|,~\|v\|=1
\}=0.$
\item [Assumption C] {\bf (Gaussian Condition)}
For any $u_i,v_i\in \mathbf h, \epsilon_i\in \mathbb Z_2: i=1,2,3$
\be \lim_{t\rightarrow 0}\frac{1}{t} \langle
\Omega,~(U_t^{(\epsilon_1)}-1)(u_1,v_1)(U_t^{(\epsilon_2)}-1)(u_2,v_2)
(U_t^{(\epsilon_3)}-1)(u_3,v_3) ~\Omega \rangle=0. \ee

\item [Assumption D] {\bf (Minimality)}\\
 The set
$ \mathcal S=\{
U_{\underbar{s},\underbar{t}}(\underbar{u},\underbar{v})
\Omega:=U_{s_1,t_1}(u_1,v_1)\cdots U_{s_n,t_n}(u_n,v_n)\Omega :
\underbar{s}=(s_1,s_2, \cdots, s_n),$ $ \underbar{t}=(t_1,t_2,
\cdots, t_n)$ $:~ 0 \le s_1\le t_1\le  s_2 \le\ldots\le s_n\le  t_n<
\infty,n\ge 1, \underbar{u}=\otimes_{i=1}^n
u_i,\underbar{v}=\otimes_{i=1}^n v_i\in \mathbf h^{\otimes n}\}$ is
total in $\mathcal H.$
\end{description}

\brmrk (a) The hypothesis {\bf A, B  and C } hold in many
situations, for  example for unitary solutions of the
Hudson-Parthasarathy flow (\ref{hpeqn})  with bounded operator
coefficients and having no Poisson  terms.\\
(b) The assumption {\bf D}  is  not really a restriction, one
 can as well work with replacing $\mathcal  H $  by span closure  of $S.$
Taking $0 \le s_1\le t_1\le s_2\le  \ldots \le s_n\le t_n< \infty$
in the definition of $S\subseteq \mathcal H$ is enough for totality
of the set $S$ because :
 for $0\le r\le s\le t \le \infty,$ we have
$U_{r,t}(p,w))=\sum_{j} U_{r,s}(p,e_j) U_{s,t}(e_j,w). $ So if there
are overlapping intervals $[s_k,t_k) \cap [s_{k+1},t_{k+1})\ne
\varnothing$ then the vector
$\xi=U_{\underbar{s},\underbar{t}}(\underbar{u},\underbar{v})
\Omega$ in $\mathcal H$ can be obtained as a vector in the closure
of the linear span of $S.$

\ermrk

\noindent For any $n\ge 1$ we have the following useful
observations.
 \blema \label{Uepbasic}
\begin{description}
\item [(i)]  For any   $  0\le r\le s\le t<\infty, $
\be U_{r,t}^{(n,k)}= U_{r,s}^{(n,k)} U_{s,t}^{(n,k)} .\ee
\item [(ii)] For any $1\le k_1,k_2,\cdots, k_m\le n~: k_i\ne k_j$ for $i\ne j$
and $ 0\le s_i\le t_i <\infty~: i=1,2, \cdots, n$ \be
\label{Uepsikn} \prod_{i=1}^{ m} U_{s_i,t_i}^{(n,\epsilon_{k_i})}
(\underbar{u},\underbar{v}) = \prod_{i=1}^{ m}
U_{s_i,t_i}^{(n,\epsilon_{k_i})} (u_{k_i},v_{k_i}) \prod_{j\ne k_i}
\langle u_j,v_j \rangle \ee for every $\underbar{u}=\otimes_{i=1}^n
u_i,\underbar{v}=\otimes_{i=1}^n v_i\in\mathbf h^{\otimes n}$ and
$\underline{\epsilon}\in \mathbb Z_2^n.$
\item [(iii)] \be
 \label{Utnbasic}
 U_{r,t}^{(n)}=U_{r,s}^{(n)}
 U_{s,t}^{(n)}.\ee
\end{description}
\elema

\begin{proof}
{\bf (i)}  It follows from the  definition and  assumptions {\bf A1} and {\bf A2}.\\
{\bf (ii)}  As in proof  of Lemma  \ref{Akn} ~{\bf (i)} by induction
(\ref{Uepsikn}) can be  proved  so  we are omitting the proof
here.\\
{\bf (iii)} Since $U_{r,t}^{(n)}$ is  a product of
$U_{r,t}^{(n,k)}~:k=1,2, \dots n$ and we have
\[U_{r,t}^{(n,k)}= U_{r,s}^{(n,k)} U_{s,t}^{(n,k)},\]
it is enough to prove that the unitary operators  $U_{r,s}^{(n,k)}$
and $ U_{s,t}^{(n,l)}$ commute  for $k\ne l.$ To see this let us
consider the following. By part  {\bf (ii)} and the fact that
$U_{r,s}(u_k,v_k) $ and $U_{s,t} (u_l,v_l)$ commute by assumption
{\bf A2}, we get
 \bean
 && U_{r,s}^{(n,k)} U_{s,t}^{(n,l)}
(\underbar{u},\underbar{v})=U_{r,s} (u_k,v_k) U_{s,t} (u_l,v_l)
\prod_{i\ne k,l} \langle u_i,v_i \rangle \\
&&=U_{s,t} (u_l,v_l)U_{r,s}(u_k,v_k) \prod_{i\ne k,l} \langle
u_i,v_i \rangle =U_{s,t}^{(n,l)} U_{r,s}^{(n,k)}
(\underbar{u},\underbar{v}) .\eean
 As
 all  the operators $U $ appear here are  bounded this implies
\[U_{r,s}^{(n,k)} U_{s,t}^{(n,l)} =U_{s,t}^{(n,l)} U_{r,s}^{(n,k)}.\]
\end{proof}

\section{Filtration}
For any  $0\le q\le  t<\infty,$  let $\mathcal H_{[q,t]}=
\overline{Span}
~\mathcal S_{[q,t]}, $ where   $ \mathcal S_{[q,t]}\subseteq \mathcal H$ is given by \\
$\{ \xi_{[q,t]}=
U_{\underbar{r},\underbar{s}}^{(n)}(\underbar{u},\underbar{v})
\Omega=U_{r_1,s_1}(u_1,v_1)\cdots U_{r_n, s_n}(u_n,v_n)\Omega \in
\mathcal S: q \le \underbar{r},\underbar{s}< t, n\ge 1,
\underbar{u},\underbar{v}\in \mathbf h^{\otimes n}\}.$  We shall
  denote the Hilbert spaces $\mathcal H_{[0,t]}$ and  $\mathcal H_{[t,\infty)}$
     by $\mathcal H_{t]}$  and $\mathcal H_{[t}$ respectively.
\blema  For $0\le  t\le  T\le \infty,$  there exist a unitary
isomorphism $\Xi_t:\mathcal H_{t]} \otimes \mathcal H_{(t,T]}
\rightarrow \mathcal H_{T]}$ such that   \be \label{Ut-filt}  U_t
(u,v)= \Xi_t^*U_t (u,v) \otimes   1_{ \mathcal H_{(t,T]}} \Xi_t.\ee
 \elema

\begin{proof}
Let us define a map $ \Xi_t :\mathcal H_{t]} \otimes \mathcal
H_{[t,T]} \rightarrow \mathcal H_{T]}$ by
\[   \Xi_t  (\xi_{[0,t]} \otimes \zeta_{[t,T]})=U_{\underbar{r},
\underbar{s}}^{(n)}(\underbar{u}, \underbar{v})
U_{\underbar{r}^\prime,\underbar{s}^\prime}^{(n)}
(\underbar{p},\underbar{w})\Omega\] for $\xi_{[0,t]}=
U_{\underbar{r},\underbar{s}}^{(n)}(\underbar{u} ,\underbar{v})
\Omega \in  \mathcal S_{t]}$ and
 $\zeta_{[t,T]}=  U_{\underbar{r}^\prime,\underbar{s}^\prime}^{(n)}
 (\underbar{p},\underbar{w})  \Omega \in  \mathcal S_{[t,T]},$  then extending linearly.\\
Now let us consider the following. By assumption {\bf A},  for
$\xi_{[0,t]}$ and $\zeta_{[t,T]}$ as above and $\eta_{[0,t]}=
U_{\underbar{a},\underbar{b}}^{(n)}(\underbar{x},\underbar{y})
\Omega \in  \mathcal S_{t]}$ and
 $\gamma_{[t,T]}=  U_{\underbar{a}^\prime,\underbar{b}^\prime}^{(n)}
 (\underbar{g},\underbar{h})  \Omega \in  \mathcal S_{[t,T]},$  we have
\bean && \langle \Xi_t (\xi_{[0,t]} \otimes \zeta_{[t,T]}), \Xi_t
(\eta_{[0,t]}
\otimes \gamma_{[t,T]}) \rangle\\
&&= \langle U_{\underbar{r},\underbar{s}}^{(n)}(\underbar{u},
\underbar{v})
U_{\underbar{r}^\prime,\underbar{s}^\prime}^{(n)}(\underbar{p},
\underbar{w})\Omega,
U_{\underbar{a},\underbar{b}}^{(n)}(\underbar{x}, \underbar{y})
U_{\underbar{a}^\prime,\underbar{b}^\prime}^{(n)}(\underbar{g},
\underbar{h}) \Omega \rangle\\
&&= \langle \Omega, \left
[U_{\underbar{r},\underbar{s}}^{(n)}(\underbar{u}, \underbar{v})
U_{\underbar{r}^\prime,\underbar{s}^\prime}^{(n)}(\underbar{p}^{(n)},\underbar{w}^{(n)})\right
]^* U_{\underbar{a},\underbar{b}}^{(n)}(\underbar{x}, \underbar{y})
U_{\underbar{a}^\prime,\underbar{b}^\prime}^{(n)}(\underbar{g},
\underbar{h}) \Omega \rangle\\
&&= \langle \Omega, \left
[U_{\underbar{r},\underbar{s}}^{(n)}(\underbar{u}, \underbar{v})
\right ]^* U_{\underbar{a},\underbar{b}}^{(n)}(\underbar{x},
\underbar{y})\Omega \rangle\\
&&~~~~ \langle \Omega, \left[
U_{\underbar{r}^\prime,\underbar{s}^\prime}^{(n)}(\underbar{p},
\underbar{w})\right ]^*
U_{\underbar{a}^\prime,\underbar{b}^\prime}^{(n)}(\underbar{g},
\underbar{h}) \Omega \rangle\\
&& =\langle \xi_{[0,t]} ,\eta_{[0,t]} \rangle \langle \zeta_{[t,T]}
,  \gamma_{[t,T]} \rangle. \eean Thus we get $\langle  \Xi_t
(\xi_{[0,t]} \otimes \zeta_{[t,T]}), \Xi_t  (\eta_{[0,t]} \otimes
\gamma_{[t,T]}) \rangle=\langle \xi_{[0,t]} \otimes \zeta_{[t,T]}~,
~ \eta_{[0,t]} \otimes \gamma_{[t,T]} \rangle.$
 Since by definition  range of  $ \Xi_t $ is dense in  $\mathcal H_{T]},$
 this  proves $ \Xi_t $ is a unitary operator.

\noindent
 Again by  similar argument as above, for any $u,v\in \mathbf h,$ we have
\bean
&& \langle  \Xi_t~ \xi_{[0,t]} \otimes \zeta_{[t,T]} ,~U_t(u,v)~ \Xi_t~ \eta_{[0,t]} \otimes \gamma_{[t,T]} \rangle\\
&&= \langle U_{\underbar{r},\underbar{s}}^{(n)}(\underbar{u},
\underbar{v}) \Omega, U_t(u,v)
U_{\underbar{a},\underbar{b}}^{(n)}(\underbar{x},
\underbar{y}) \Omega \rangle\\
&&~~~~ \langle
U_{\underbar{r}^\prime,\underbar{s}^\prime}^{(n)}(\underbar{p},
\underbar{w}) \Omega ,
U_{\underbar{a}^\prime,\underbar{b}^\prime}^{(n)}(\underbar{g},
\underbar{h}) \Omega \rangle\\
&&= \langle \xi_{[0,t]} ~,~U_t(u,v) \eta_{[0,t]} \rangle~~ \langle
\zeta_{[t,T]} ,\gamma_{[t,T]} \rangle \eean This  proves ~
(\ref{Ut-filt}).
\end{proof}

\section{Expectation semigroups}
Let us look at the various  semigroups  associated with the unitary
evolution $\{U_{s,t}\}.$

\noindent For any fixed  $n\ge 1,$ we define a family of  operators
$\{T_t^{(n)}\}$ on $\mathbf h ^{\otimes n}$ by setting
\[\langle \phi,T_t^{(n)} ~\psi \rangle:= \langle \Omega,
  U_t^{(n)}(\phi,\psi)~ \Omega \rangle,~ \forall
  \phi,\psi \in { \mathbf h}^{\otimes n}.\]
Then in particular  for product vectors $\underbar{u}=
\otimes_{i=1}^n u_i,~ \underbar{v}=\otimes_{i=1}^n v_i\in\mathbf h
^{\otimes n}$
\[\langle \underbar{u},T_t^{(n)} ~\underbar{v}   \rangle= \langle \Omega,
U_t^{(n)}(\underbar{u},\underbar{v})~ \Omega \rangle = \langle
\Omega,   U_t(u_1,v_1) U_t(u_2,v_2)\cdots U_t(u_n,v_n)~ \Omega
\rangle.\] For $n=1,$  we shall write $T_t$  for the family
$T_t^{(1)}.$

\blema
 The above family of operators $\{T_t^{(n)}\}$ is  a   semigroup of contractions on
   $\mathbf h ^{\otimes n}.$
\elema

\begin{proof} Since $U_t^{(n)}$ is in particular  contractive,   for any
$ \phi,\psi \in
 { \mathbf h}^{\otimes n}$
\[|\langle \phi,T_t^{(n)} ~\psi \rangle|= |\langle \phi~\Omega,
   U_t^{(n)}  \psi~ \Omega \rangle| \le\|\phi\|
   ~\|\psi \|\]  and contractivity of  $T_t^{(n)}$ follows.

\noindent In order to prove that   this family of contractions
$T_t^{(n)}$ is a semigroup it is enough to show that   for any  $0
\le s \le t $ and product vectors $\underbar{u}= \otimes_{i=1}^n
u_i,~ \underbar{v}=\otimes_{i=1}^n v_i\in\mathbf h ^{\otimes n},$
\[\langle \underbar{u} ,T_t^{(n)} ~\underbar{v}   \rangle=\langle
 \underbar{u}, T_s^{(n)} T_{t-s}^{(n)}  \underbar{v}\rangle.\]
Consider the  product orthonormal basis   $\{
e_{\underbar{j}}=e_{j_1} \otimes e_{j_2} \otimes \cdots \otimes
e_{j_n} :~\underbar{j}=(j_1,j_2,\cdots, j_n): j_1,j_2,\cdots, j_n\ge
1 \}$ of $\mathbf h ^{\otimes n}.$ By part (iii) of Lemma \ref{Auv}
and evolution property (\ref{Utnbasic}) of $U_t^{(n)}$ ,
 \bean
&&\langle \underbar{u} ,T_t^{(n)} ~\underbar{v}  \rangle
=\langle \Omega,   U_t^{(n)}(\underbar{u},\underbar{v})~ \Omega \rangle\\
&&=\sum_{ \underbar{j} }
\langle   \Omega,~ U_s^{(n)}(\underbar{u} ,e_{\underbar{j}} ) U_{s,t}^{(n)}(e_{\underbar{j}},  \underbar{v}) \Omega\rangle\\
&&=\sum_{ \underbar{j} } \langle   \Omega,~
U_s^{(n)}(\underbar{u},e_{\underbar{j}})~\Omega \rangle
\langle  \Omega~,   U_{t-s}^{(n)}(e_{\underbar{j}},  \underbar{v}) \Omega\rangle\\
&&=\sum_{ \underbar{j} } \langle \underbar{u}, T_s^{(n)}
e_{\underbar{j}} \rangle \langle e_{\underbar{j}},    T_{t-s}^{(n)}
\underbar{v}\rangle=\langle   \underbar{u} , T_s^{(n)} T_{t-s}^{(n)}
\underbar{v}\rangle. \eean

\end{proof}

\noindent The following Lemma will  be needed  in the sequel \blema
\label{UTkn}
\begin{description}
\item[(i)] For $1\le k\le n,$
 \be
\langle \Omega ,U_t^{(n,k)} (\underbar{p}, \underbar{w}) \Omega
\rangle =\langle \underbar{p}, 1_{{\mathbf h} ^{(\otimes k-1)}}
\otimes T_t \otimes 1_{{\mathbf h}^{(\otimes n-k)}} \underbar{w}
\rangle, ~\forall \underbar{p},\underbar{w} \in { \mathbf
h}^{\otimes n}. \ee We shall denote the ampliation $ 1_{{\mathbf
h}^{(\otimes k-1)}} \otimes T_t \otimes 1_{{\mathbf h}^{(\otimes
n-k)}}$  by $T_t^{(n,k)}  $.
\item[(ii)] For  any $1\le m\le n,~ \underbar{p},\underbar{w}
\in { \mathbf h}^{\otimes n},$
\[
 \langle \Omega ,(\prod_{k=1}^m U_t^{(n,k)})  (\underbar{p}, \underbar{w}) \Omega \rangle
=\langle  \underbar{p}, T_t^{(m)}\otimes 1_{{\mathbf h}^{(\otimes
n-m)} }~ \underbar{w} \rangle. \]
\item[(iii)]
 For any $\phi\in { \mathbf h}^{\otimes n},$
\bean
&&\|(U_t^{(n,k)}-1) \phi\otimes \Omega \|^2\\
&& = \langle (1-T_t^{(n,k)}) \phi,\phi \rangle
+\langle \phi,(1-T_t^{(n,k)}) \phi\rangle \\
&&\le 2\|1-T_t\|~\|\phi\|^2. \eean
\item[(iv)]   \bean
&&
\|(U_t^{(n)} -1) \phi\otimes \Omega \|^2\\
&& = \langle (1-T_t^{(n)}) \phi,\phi\rangle
+\langle \phi,(1-T_t^{(n)}) \phi\rangle \\
&&\le 2\|(1-T_t^{(n)})\|~\|\phi\|^2.\eean
\item[(v)] For any $v\in \mathbf h$
\be \sum_{m\ge 1} \|(U_t-1)(e_m,v)\Omega\|^2=2 Re \langle v ,
(1-T_t) v \rangle \le 2 \|v\|^2  \|T_t-1\|.\ee

\end{description}
\elema
\begin{proof}

(i) It follows from  the fact that for product vectors \be \langle
\Omega ,U_t^{(n,k)} (\underbar{p},\underbar{w}) \Omega \rangle =
\langle p_k, T_t^{(n,k)} w_k \rangle \prod_{i\ne k}\langle  p_i, w_i
\rangle . \ee The part {\bf (ii)} follows from Lemma \ref{Uepbasic}
{\bf (ii)}.

 \noindent Proof of {\bf (iii)} and  {\bf (iv)} are similar  so we
prove only for $U_t^{(n,k)}.$ We have \bean
&&  \|(U_t^{(n,k)} -1)~\phi\Omega \|^2\\
&&= \langle \phi\Omega,   [(U_t^{(n,k)} -1)^*(U_t^{(n,k)}
-1)]\phi\Omega\rangle
\\
&&\le \langle \phi  \Omega, [2-(U_t^{(n,k)})^* -U_t^{(n,k)}]
\phi\Omega\rangle~~~\mbox{(since $U_t^{(n,k)}$ is in particular
contractive)}
\\
&& = \langle (1-T_t^{(n,k)}) \phi,\phi \rangle +\langle
\phi,(1-T_t^{(n,k)})\phi\rangle
\eean Thus the statement  follows.\\

{\bf (v)}  For any $v\in \mathbf h$
 \bean &&\sum_{m\ge 1}
\|(U_t-1)(e_m,v)\Omega\|^2\\
&&= \sum_{m\ge 1} \langle \Omega
,(U_t-1)^*(v,e_m))(U_t-1)(e_m,v)\Omega \rangle\\
&&=\langle \Omega ,[(U_t-1)^*(U_t-1)](v,v)\Omega \rangle\\
&&\le\langle \Omega ,[2-U_t^*-U_t](v,v)\Omega \rangle\\
&&=\langle v ,[2-T_t^*-T_t]v \rangle=2 Re \langle v , (1-T_t) v
\rangle \le 2 \|v\|^2  \|T_t-1\|. \eean
\end{proof}
\noindent Now we are ready to prove \bpro
  Under the  assumption {\bf{B}} the semigroup  $\{T_t^{(n)} \}$ is  uniformly continuous.
\epro

\begin{proof}
\noindent Assumption   {\bf B } on the family of unitary operators
$\{U_{s,t}\}$ implies that the semigroup of contractions $ \{T_t\}$
on $\mathbf h$ is   uniformly continuous. To apply  induction let us
assume that for some $m\ge 1,$ the  contractive semigroups
$\{T_t^{(n)}\} $ are uniformly continuous for all $1\le n \le m-1.$
Now, for any $ \phi,\psi\in {\mathbf h}^{\otimes m}$ \bean
&& \langle \phi \otimes \Omega,   (U_t^{(m)} -1) \psi \otimes \Omega \rangle\\
&& =\langle \phi\otimes \Omega,   \left (
[\prod_{k=1}^{m-1}U_t^{(m,k)}]
 [U_t^{(m,m)}]-1 \right )\psi\otimes  \Omega \rangle\\
&& =\langle  [\prod_{k=1}^{m-1}U_t^{(m,k)}] ^*\phi\otimes \Omega,
\left ([U_t^{(m,m)}]-1 \right )\psi\otimes  \Omega \rangle\\
&& ~~~+\langle \phi\otimes \Omega,   \left (
[\prod_{k=1}^{m-1}U_t^{(m,k)}]-1 \right )\psi\otimes  \Omega
\rangle.\eean Taking absolute value, by Lemma \ref{UTkn}~we get
\bean
&& |\langle \phi,~ (T_t^{(m)}-1_{{\mathbf h}^{\otimes m }})\psi \rangle|\\
&& \le \| \phi\|~ \|\psi\|    \sqrt{2\|T_t^{(m,m)}- 1_{{\mathbf
h}^{\otimes m }} \|}
 + |\langle \phi,   \left ( [T_t^{(m-1)} \otimes 1_{\mathbf h}]-
 1_{{\mathbf h}^{\otimes m }} \right )\psi \rangle|\\
&&\le \|\phi\|  \|\psi\| \left [\sqrt{2\|T_t-1\|} + \| T_t^{(m-1)}-1
\| \right ]. \eean So uniform continuity of $T_t^{(m-1)}$ and $T_t$
implies that $T_t^{(m)}$ is uniformly continuous.
 \end{proof}

\noindent Let us denote the bounded generator of the uniformly
continuous semigroup $T_t^{(n)}$ on $\mathbf h ^{\otimes n}$  by
$G^{(n)}$ and for $n=1$ by  $G.$

\noindent For $m,n\ge 1,$ we define a family of operators $\{Z_t^{(
m,n)}:t\ge 0\}$ on the Banach space $\mathcal B_1( \mathbf
h^{\otimes m}, \mathbf h^{\otimes n})$ by

\[Z_t^{(   m,n)} \rho=  Tr_{\mathcal H}[ U_t^{(n)} ( \rho \otimes |\Omega ><  \Omega|) (U_t^{(m)})^*],~\rho \in \mathcal B_1( \mathbf h^{\otimes m}, \mathbf h^{\otimes n}) .\]
Then in particular for product vectors $\underbar{u},\underbar{v}
\in \mathbf h^{\otimes m},\underbar{p},\underbar{w} \in \mathbf
h^{\otimes n}.$

\be \label{Ztsimp}\langle \underbar{p} ,Z_t^{(   m,n)}
(|\underbar{w}><\underbar{v}|) \underbar{u}  \rangle:= \langle
 U_t^{(m)}(\underbar{u},\underbar{v}) \Omega,~ U_t^{(n)}(\underbar{p},\underbar{w})~
\Omega \rangle.\ee

\blema The above family  $\{Z_t^{( m,n)}\}$ is a  semigroup of
contractive maps on $\mathcal B_1( \mathbf h^{\otimes m}, \mathbf
h^{\otimes n}).$  Furthermore, assumption {\bf {B}} implies that
$\{Z_t^{( m,n)}\}$ is uniformly continuous. \elema
\begin{proof}
For $\rho\in \mathcal B_1( \mathbf h^{\otimes m}, \mathbf h^{\otimes
n})$
 \bean
&&\|Z_t^{(m,n)} \rho\|_1=  \|Tr_{\mathcal H}[ U_t^{(n)} ( \rho \otimes |\Omega ><  \Omega|) (U_t^{(m)})^*]\|_1\\
&& = \sup\{\sum_{k\ge 1} |\langle \phi_k^{(n)} ,Tr_{\mathcal H}[
U_t^{(n)}( \rho \otimes |\Omega ><  \Omega|) (U_t^{(m)})^*]
\phi_k^{(m)} \rangle|~:~\{\phi_k^{(l)}\}\\
&&~~~~~~~~~~~~~~~~~~~~~~~~~~~~~~~~ is ~an~ ONB~ of~
\mathbf h ^{\otimes l}, l=m,n \}\\
&& \le  \sup_{\phi^{(l)}} \sum_{j,k\ge 1} |\langle
\phi_k^{(n)}\otimes \zeta_j,
 U_t^{(n)} ( \rho \otimes |\Omega ><  \Omega|) (U_t^{(m)})^* \phi_k^{(m)}
 \otimes \zeta_j \rangle|\\
&& \le \|U_t^{(n)} ( \rho \otimes |\Omega >< \Omega|)
(U_t^{(m)})^*\|_1. \eean Since for any $l\ge 1,~\{U_t^{(l)}\}$ is in
particular a contractive  family of operators
\[\|Z_t^{(m,n)} \rho\|_1 \le  \|\rho \otimes |\Omega ><  \Omega|\|_1 =  \|\rho\|_1.\]

\noindent In order to prove that the family of contractions
$\{Z_t^{(m,n)}\}$  is a semigroup it is enough to verify that
property  for the rank one  operator
$\rho=|\underbar{w}><\underbar{v}| ~:~ \underbar{w}=\otimes_{i=1}^n
w_i\in { \mathbf h}^{\otimes n},~ \underbar{v}=\otimes_{i=1}^m v_i
\in {\mathbf h}^{\otimes m}.$  Therefore, it suffices to prove that
for $\underbar{p}=\otimes_{i=1}^n p_i\in { \mathbf h}^{\otimes n},~
\underbar{u}=\otimes_{i=1}^m u_i\in { \mathbf h}^{\otimes m}$
\[\langle \underbar{p},   Z_t^{(m,n)} (\rho)    \underbar{u} \rangle=
\langle \underbar{p}, Z^{(m,n)}_s Z^{(m,n)}_{t-s} (\rho)
\underbar{u} \rangle .\] By Lemma \ref{Uepbasic},  part (iv) of
Lemma \ref{Auv} and assumption {\bf{A}} for $0 \le s \le t ,~
\underbar{u},\underbar{v} \in \mathbf h^{\otimes m}$ and product ONB
$\{ e_{\underbar{j}}^{(m)}=e_{j_1}\otimes \cdots \otimes e_{j_m} \}$
of $\mathbf h^{\otimes m}$ and $\{ e_{\underbar{k}}^{(n)}=e_{k_1}
\otimes \cdots \otimes e_{k_n} \}$ of $\mathbf h^{\otimes n}$

\bean
&& \langle   U_t^{(m)}(\underbar{u},\underbar{v}) \Omega,~ U_t^{(n)}(\underbar{p},\underbar{w})~\Omega \rangle\\
&&=\sum_{j,k} \langle U_s^{(m)}(\underbar{u},e_{\underbar{j}}^{(m)})
\Omega,~U_s^{(n)} (\underbar{p},e_{\underbar{k}}^{(n)} )~\Omega
\rangle \langle U_{t-s}^{(m)}(e_{\underbar{j}}^{(m)}, \underbar{v})
\Omega,~U_{t-s}^{(n)}(e_{\underbar{k}}^{(n)}, \underbar{p} )~\Omega
\rangle. \eean This give \bean
&& \langle \underbar{p} ,Z_t^{(m,n)} (\rho) \underbar{u}  \rangle\\
 &&=\sum_{\underbar{j},\underbar{k}}
\langle \underbar{p}, Z^{(m,n)}_s (| e_{\underbar{k}}^{(n)} ><e_{\underbar{j}}^{(m)} |)   \underbar{u} \rangle \langle  e_{\underbar{k}}^{(n)} ,  Z^{(m,n)}_{t-s} (\rho) e_{\underbar{j}}^{(m)}\rangle\\
&&=\sum_{\underbar{j},\underbar{k}}
\langle   e_{\underbar{j}}^{(m)},  (Z^{(m,n)}_s)^*   (|\underbar{u}><\underbar{p}|) e_{\underbar{k}}^{(n)}   \rangle \langle  e_{\underbar{k}}^{(n)} ,  Z^{(m,n)}_{t-s} (\rho) e_{\underbar{j}}^{(m)} \rangle\\
&&=\sum_{\underbar{j}}
\langle   e_{\underbar{j}}^{(m)},  (Z^{(m,n)}_s)^*   (|\underbar{u}><\underbar{p}|) Z^{(m,n)}_{t-s}(\rho) e_{\underbar{j}}^{(m)}\rangle\\
&&=Tr [(Z^{(m,n)}_s)^*   (|\underbar{u}><\underbar{p}|) Z^{(m,n)}_{t-s}(\rho)]\\
&&=Tr [ |\underbar{u}><\underbar{p}|   Z^{(m,n)}_s Z^{(m,n)}_{t-s} (\rho)]\\
&&=
\langle \underbar{p},   Z^{(m,n)}_s Z^{(m,n)}_{t-s} (\rho)   \underbar{u} \rangle. \\
\eean

\noindent In order to prove uniform  continuity of  $ Z_t^{( m,n)}$
we  consider
 \bean
&&\|(Z_t^{(   m,n)}-1) (|\underbar{w}><\underbar{v}|)\|_1\\
&& = \sup\{\sum_{k\ge 1} |\langle \phi_k^{(n)} ,
(Z_t^{(   m,n)} -1) (|\underbar{w}><\underbar{v}|)  \phi_k^{(m)} \rangle|: \{\phi_k^{(l)}\}\\
&&~~~~~~~~~~~~~~~~~~~~~~~~~~~~~~~~~~~is ~an~ ONB~ of~
\mathbf h ^{\otimes l}, l=m,n \}\\
&& = \sup_{\phi^{(l)}}
 \sum_{k\ge 1} |\langle U_t^{(m)} ( \phi_k^{(m)},\underbar{v}) \Omega    ,
 U_t^{(n)}  ( \phi_k^{(n)}, \underbar{w} ) \Omega  \rangle-\overline{\langle
  \phi_k^{(m)},\underbar{v} \rangle} \langle \phi_k^{(n)}, \underbar{w}   \rangle|\\
&& \le \sup_{\phi^{(l)}} \sum_{k\ge 1} |\langle (U_t^{(m)}-1)
 ( \phi_k^{(m)},\underbar{v}) \Omega    , U_t^{(n)}  ( \phi_k^{(n)}, \underbar{w} )
 \Omega  \rangle|\\
&& +\sup_{\phi^{(l)}}  \sum_{k\ge 1} |\overline{\langle
\phi_k^{(m)},
\underbar{v} \rangle}  \langle  \Omega    , (U_t^{(n)}-1)  ( \phi_k^{(n)}, \underbar{w} ) \Omega  \rangle|\\
&& \le \sup_{\phi^{(l)}} \left [\sum_{k\ge 1} \| (U_t^{(m)}-1) ( \phi_k^{(m)},\underbar{v})
\Omega \|^2 \right ]^{\frac{1}{2}}
 \left [\sum_{k\ge 1} \|U_t^{(n)}  ( \phi_k^{(n)}, \underbar{w} ) \Omega \|^2 \right ]^{\frac{1}{2}} \\
&& +\sup_{\phi^{(l)}} \left [\sum_{k\ge 1} |\langle
\phi_k^{(m)},\underbar{v} \rangle|^2 \right ]^{\frac{1}{2}}
\left [\sum_{k\ge 1} \|(U_t^{(n)}-1)  ( \phi_k^{(n)}, \underbar{w} ) \Omega \|^2 \right ]^{\frac{1}{2}} \\
\eean
 So by Lemma \ref{UTkn}
\bean
&&\|(Z_t^{(   m,n)}-1) (|\underbar{w}><\underbar{v}|)\|_1\\
&&\le \sqrt{2} \|\underbar{v}\| ~ \|\underbar{w}\|
\left(\sqrt{\|T_t^{(m)}-1\|}+\sqrt{\|T_t^{(n)}-1\|}\right) \eean

Now for any  $\rho=\sum_{k} \lambda_k |
\phi_k^{(n)}><\phi_k^{(\underline{\alpha})}|\in \mathcal B_1(
\mathbf h^{\otimes m}, \mathbf h^{\otimes n})$  we have \bean
&&\|Z_t^{(   m,n)} (\rho)- \rho\|_1\\
&&\le \sqrt{2} \sum_{k} | \lambda_k|
\left(\sqrt{\|T_t^{(m)}-1\|}+\sqrt{\|T_t^{(n)}-1\|}\right)
\\
&&\le \sqrt{2}
\|\rho\|_1~\left(\sqrt{\|T_t^{(m)}-1\|}+\sqrt{\|T_t^{(n)}-1\|}\right).
\eean
\end{proof}

Thus by  uniform  continuity of the semigroup $ T_t^{(m)}$ and
$T_t^{(n)}$ it follows that  the semigroup ${Z_t^{( m,n)}}$ is
uniformly continuous on $\mathcal B_1( \mathbf h^{\otimes m},
\mathbf h^{\otimes n}).$

We shall  denote the  bounded  generator of  the semi-group $Z_t^{(
m,n)}$ by $\mathcal L^{( m,n)}.$ For $n\ge 1$ we shall write
$Z_t^{(n)}$ for   the semi-group $Z_t^{( n,n)}$ on the Banach space
$\mathcal B_1(\mathbf h^{\otimes n})$ and shall  denote  its
generator by $\mathcal L^{(n)}.$ Moreover,  we denote the semigroup
$Z_t^{(1)}$ and its generator  $\mathcal L^{(1)}$ by just $Z_t$  and
$\mathcal L$  respectively.

\blema  \label{Zn+} For any $n\ge 1 ,~Z_t^{(n)}$ is a positive trace
preserving   semigroup on $\mathcal B_1(\mathbf h ^{\otimes n}).$
 \elema
\begin{proof}
Positivity follows from
\[ \langle \underbar{u} ,Z_t^{(n)} (|\underbar{v}><\underbar{v}|) \underbar{u}  \rangle = \|  U_t^{(n)}(\underbar{u},\underbar{v}) \Omega\|^2 \ge 0~
\forall~~ \underbar{u},\underbar{v} \in \mathbf h ^{\otimes n}.\] To
prove that  $Z_t^{(n)}$ is trace preserving it is enough to show
that
\[Tr [Z_t^{(n)} (|\underbar{u}><\underbar{v}|)]= \langle  \underbar{v},  \underbar{u} \rangle. \]
By definition and  Lemma \ref{Auv}  \bean
&& Tr[ Z_t^{(n)} (|\underbar{u}><\underbar{v}|)]= \sum_{k}\langle \underbar{e}_k ,Z_t^{(n)} (|\underbar{u}><\underbar{v}|) \underbar{e}_k  \rangle\\
&&=\sum_{k}\langle U_t^{(n)}(\underbar{e}_k,\underbar{v})\Omega ,U_t^{(n)}(\underbar{e}_k,\underbar{u}) \Omega \rangle\\
&&=\langle \Omega ,(U_t^{(n)})^*
U_t^{(n)}(\underbar{v},\underbar{u}) \Omega \rangle. \eean Since
$U_t^{(n)}$ is   unitary, we get
\[Tr[ Z_t^{(n)} (|\underbar{u}><\underbar{v}|)]=\langle  \underbar{v},  \underbar{u} \rangle.\]
\end{proof}

\noindent This Lemma  gives \be \label{tr-zero} Tr(\mathcal
L^{(n)}\rho)=0,~\forall  \rho\in  \mathcal B_1(\mathbf h^{\otimes
n}). \ee

\noindent  We also need another class of semigroup. For $m,n\ge 1$
we define a family of maps $F^{(m,n)}_t$ on the Banach space
$\mathcal B_1( \mathbf h^{\otimes m}, \mathbf h^{\otimes n})$ by
 \be F^{(m,n)}_t \rho=  Tr_{\mathcal H}[
(U_t^{(n)})^*
 ( \rho \otimes |\Omega ><  \Omega|) U_t^{(m)}], ~\forall  \rho
 \in  \mathcal B_1( \mathbf h^{\otimes m}, \mathbf h^{\otimes n})
 \ee

\noindent So in particular for product vectors
$\underbar{u},\underbar{v} \in \mathbf h^{\otimes
m},\underbar{p},\underbar{w} \in \mathbf
h^{\otimes n},$\\
  $\langle
\underbar{p} ,F^{(m,n)}_t (|\underbar{w}><\underbar{v}|)
\underbar{u}
  \rangle
=\langle   ( U_t^{(m)})^*(\underbar{u},\underbar{v}) \Omega,~ (
U_t^{(n)})^* (\underbar{p},\underbar{w})~\Omega \rangle. $ \blema
For any $m,n\ge 1, ~\{F^{(m,n)}_t:t\ge 0\} $ is a uniformly
continuous contractive semigroup. \elema
\begin{proof} Similarly as for the semigroup  ${Z_t^{( m,n)}}.$
\end{proof}
For $n\ge 1,$ we shall write $F_t^{(n)}$ for   the semi-group
$F_t^{( n,n)}$ on the Banach space $\mathcal B_1(\mathbf h^{\otimes
n})$ and in particular  $F_t$  for the semigroup  $F_t^{(1)}$ on
$\mathcal B_1(\mathbf h).$ We conclude this section by the
following useful observation.

\blema Under the {\bf Assumption C}, for any  $n\ge
3,~\underbar{u},\underbar{v} \in \mathbf h^{\otimes
n},\underline{\epsilon}\in \mathbb Z_2^n$
 \be \label{4Ut-1} \lim_{t\rightarrow
0}\frac{1}{t} \langle \Omega,~ (U_t^{(\epsilon_1)}-1)(u_1,v_1)\cdots
(U_t^{(\epsilon_n)}-1)(u_n,v_n)~\Omega \rangle=0 .\ee \elema
\begin{proof}
 We have
 \bean
  &&|\frac{1}{t}\langle
[(U_t^{(\epsilon_1)}-1)(u_1,v_1)(U_t^{(\epsilon_2)}-1)(u_2,v_2)]^*~\Omega,~\\
&&~~~~~~~~~~~~~~~~ (U_t^{(\epsilon_1)}-3)(u_3,v_3)\cdots
(U_t^{(\epsilon_n)}-1)(u_n,v_n)~\Omega \rangle|^2\\
&&\le \frac{1}{t^2} \|
[(U_t^{(\epsilon_1)}-1)(u_1,v_1)(U_t^{(\epsilon_2)}-1)(u_2,v_2)]^*~\Omega\|^2\\
&&~~~~~~~~\|(U_t^{\epsilon_3}-1)(u_3,v_3)\cdots
(U_t^{(\epsilon_n)}-1)(u_n,v_n)
~\Omega\|^2\\
&&\le  C_{\underbar{u},\underbar{v}} ~~\frac{1}{t} \|
[(U_t^{(\epsilon_1)}-1)(u_1,v_1)(U_t^{(\epsilon_2)}-1)(u_2,v_2)]^*~\Omega\|^2\\
&&~~~~~~~~\frac{1}{t} \|(U_t^{(\epsilon_{n-1})}-1)(u_{n-1},v_{n-1})
(U_t^{(\epsilon_n)}-1)(u_n,v_n) ~\Omega\|^2 \eean
 for some constant $C_{\underbar{u},\underbar{v}}$ independent of  $t.$ So to prove (\ref{4Ut-1}) it is enough to show that
 for any $u,v,p,w \in \mathbf h$  and $\epsilon, \epsilon^\prime \in \mathbb Z_2$
 \be \label{UtUt*}\lim_{t\rightarrow 0}\frac{1}{t} \|
(U_t^{(\epsilon)}-1)(u,v)(U_t^{(\epsilon^\prime)}-1)(p,w)~\Omega\|^2=0.\ee
So let us look at the following
 \bean &&\|
(U_t^{(\epsilon)}-1)(u,v)(U_t^{(\epsilon^\prime)}-1)(p,w)~\Omega\|^2\\
&&= \langle
(U_t^{(\epsilon)}-1)(u,v)(U_t^{(\epsilon^\prime)}-1)(p,w)~\Omega,
(U_t^{(\epsilon)}-1)(u,v)(U_t^{(\epsilon^\prime)}-1)(p,w)~\Omega \rangle \\
&&= \langle (U_t^{(\epsilon^\prime)}-1)(p,w)~\Omega,
[(U_t^{(\epsilon)}-1)(u,v)]^*
(U_t^{(\epsilon)}-1)(u,v)(U_t^{(\epsilon^\prime)}-1)(p,w)~\Omega
\rangle. \eean \noindent By part {\bf (v)}  of Lemma {\ref {Auv}}
the above quantity is \\
$\le \|u\|^2 \langle (U_t^{(\epsilon^\prime)}-1)(p,w)~\Omega,
[(U_t^{(\epsilon)}-1)^*
(U_t^{(\epsilon)}-1)](v,v)(U_t^{(\epsilon^\prime)}-1)(p,w)~\Omega
\rangle.$ Since  by contractivity of the family
$U_t^{(\epsilon)},~~(U_t^{(\epsilon)})^*U_t^{(\epsilon)}\le 1,$ we
get \bean &&\|
(U_t^{(\epsilon)}-1)(u,v)(U_t^{(\epsilon^\prime)}-1)(p,w)~\Omega\|^2\\
 &&\le \|u\|^2 \langle
(U_t^{(\epsilon^\prime)}-1)(p,w)~\Omega,
[1-(U_t^{(\epsilon)})^*+1-U_t^{(\epsilon)}] (v,v)(U_t^{(\epsilon^\prime)}-1)(p,w)~\Omega \rangle \\
&&=- \|u\|^2 \langle (U_t^{(\epsilon^\prime)}-1)(p,w)~\Omega,
[U_t^{(1+\epsilon)}-1] (v,v)(U_t^{(\epsilon^\prime)}-1)(p,w)~\Omega \rangle \\
&&~~~~- \|u\|^2 \langle (U_t^{(\epsilon^\prime)}-1)(p,w)~\Omega,
(U_t^{(\epsilon)}-1) (v,v)(U_t^{(\epsilon^\prime)}-1)(p,w)~\Omega
\rangle. \eean Thus by {\bf Assumption C} we get (\ref{UtUt*}) and
the proof is complete.
 \end{proof}
\section{Representation of Hilbert tensor algebra and Hudson-Parthasarathy (HP) equation }
\noindent We define a scalar valued  map $K$ on $M\times M $ by
setting, for $(\underbar{u} ,\underbar{v},\underline{\epsilon}),
(\underbar{p},\underbar{w},\underline{\epsilon}^\prime) \in M_0,$
\[K\left( (\underbar{u},\underbar{v},\underline{\epsilon}), (\underbar{p},
\underbar{w},\underline{\epsilon}^\prime)   \right
):=\lim_{t\rightarrow 0}\frac{1}{t} \langle
(U_t^{(\underline{\epsilon})}-1)(\underbar{u}, \underbar{v})
\Omega,~ (U_t^{\underline{\epsilon}^\prime}-1)(\underbar{p},
\underbar{w})~\Omega \rangle,~\mbox{when  it exists.}~\]

\blema \label{kolmo} {\bf(i)}The map $K$ is a well defined  positive
definite
kernel on $M.$\\
{\bf(ii)} Up to unitary equivalence there exists a
unique separable Hilbert space $\mathbf k,$ an  embedding $\eta:
M\rightarrow \mathbf k$ and a $*$-representation $\pi$ of  $M,~
\pi:M\rightarrow \mathcal B( \mathbf k)$ such that
 \be \label{eta-dense}
 \{\eta (\underbar{u},\underbar{v},\underline{\epsilon}): (\underbar{u},
 \underbar{v},\underline{\epsilon})\in M_0\} ~\mbox{ is total in}~
 \mathbf k,\ee \be \label{eta-kelnel-} \langle \eta
(\underbar{u},\underbar{v},\underline{\epsilon}),
\eta(\underbar{p},\underbar{w},\underline{\epsilon}^\prime)  \rangle
= K\left( (\underbar{u},\underbar{v},\underline{\epsilon}),
(\underbar{p},\underbar{w},\underline{\epsilon}^\prime)   \right
)\ee and
 \be \label{repn} \pi
(\underbar{u},\underbar{v},\underline{\epsilon})\eta
(\underbar{p},\underbar{w},\underline{\epsilon}^\prime)=\eta
(\underbar{u}\otimes \underbar{p},\underbar{v}\otimes
\underbar{w},\underline{\epsilon}\oplus
\underline{\epsilon}^\prime)- \langle
\underbar{p},\underbar{w}\rangle \eta
(\underbar{u},\underbar{v},\underline{\epsilon}).\ee
 \elema
\begin{proof} {\bf(i)}
First note that for any
$(\underbar{u},\underbar{v},\underline{\epsilon})\in M_0,$
$\underbar{u}=\otimes_{i=1}^n u_i,\underbar{v}=\otimes_{i=1}^n
v_i,\\
\underline{\epsilon}=(\epsilon_1,\epsilon_2,\cdots, \epsilon_n)$ we
can write
 \bea \label{Ut-1-decompose}
&&(U_t^{(\underline{\epsilon})}-1)(\underbar{u},\underbar{v})=
\prod_{i=1}^n U_t^{(\epsilon_i)}(u_i,v_i)-\prod_{i=1}^n \langle u_i,v_i \rangle {\nonumber}\\
&&=\sum_{1\le i\le n} (U_t-1)^{(\epsilon_i)}(u_i,v_i)\prod_{j\ne i} \langle u_j,v_j \rangle{\nonumber}\\
&&+\sum_{2\le l\le n} ~~\sum_{1\le i_1<\ldots < i_m\le n}
\prod_{k=1}^l(U_t-1)^{\epsilon_{i_k}}(u_{i_k},v_{i_k}) \prod_{j\ne
i_k} \langle u_j,v_j \rangle.\eea

\noindent Now by Lemma \ref{4Ut-1}, for elements
$(\underbar{u},\underbar{v},\underline{\epsilon}),
(\underbar{p},\underbar{w},\underline{\epsilon}^\prime) \in M_0,$
$\underline{\epsilon}\in \mathbb Z_2^m$  and
$\underline{\epsilon}^\prime\in \mathbb Z_2^n,$
 we have
\bean && K\left( (\underbar{u},\underbar{v},\underline{\epsilon}),
(\underbar{p}, \underbar{w},\underline{\epsilon}^\prime)   \right
)=\lim_{t\rightarrow 0}\frac{1}{t} \langle
(U_t^{(\underline{\epsilon})}-1)(\underbar{u}, \underbar{v})
\Omega,~ (U_t^{\underline{\epsilon}^\prime}-1)(\underbar{p},
\underbar{w})~\Omega \rangle\\
&&= \sum_{1\le i\le m,~1\le j\le n}
 \prod_{k\ne i} \overline{\langle u_k,v_k \rangle}
 \prod_{l\ne j} \langle p_l,w_l \rangle
 \lim_{t\rightarrow 0}\frac{1}{t}
\langle (U_t-1)^{(\epsilon_i)}(u_i,v_i) ~\Omega,
(U_t-1)^{\epsilon_j^\prime}(p_j,w_j) ~\Omega\rangle. \eean Hence
existence of the above  limit  follows from the fact that the
semigroups $T_t$ on $\mathbf h$  and $Z_t, F_t$ on $\B_1(\mathbf h)$
are uniformly  continuous  and
 \bean && \langle (U_t-1)^{(\epsilon_i)}(u_i,v_i) ~\Omega,
(U_t-1)^{\epsilon_j^\prime}(p_j,w_j) ~\Omega\rangle\\
&&=\langle U_t^{(\epsilon_i)}(u_i,v_i) \Omega,~
U_t^{\epsilon_j^\prime}(p_j,w_j) ~\Omega \rangle- \overline{\langle
u_i,v_i\rangle} \langle
p_j,w_j \rangle \\
&&- \overline{\langle  u_i,v_i \rangle} \langle    \Omega,~
 [(U_t^{\epsilon_j^\prime}-1)(p_j,w_j)]~
 \Omega \rangle\\
&&-\overline{\langle \Omega,~ [(U_t^{(\epsilon_i)}-1)(u_i,v_i)]
 \Omega \rangle} \langle   p_j,w_j\rangle.
\eean Thus $K$ is well defined on $M_0.$ Now  extend this to the
$*$-algebra $M$ sesqui-linearly.\\
In particular we have
 \bea \label{kernel}
&& K( (u,v,0), (p,w,0) ) {\nonumber} \\
&&=\lim_{t\rightarrow 0}\{ \langle p,\frac{Z_t-1}{t}(|w><v|)u\rangle
 - \overline{\langle u, v \rangle }~\langle  p,\frac{T_t-1}{t} w \rangle
-\overline{\langle u, \frac{T_t-1}{t}v \rangle}~\langle   p,w
\rangle \} {\nonumber}\\
 &&= \langle p,\mathcal L (|w><v|) u \rangle-\overline{\langle u,v
\rangle} \langle p,G~w  \rangle- \overline{\langle u, G ~v
\rangle}\langle p,w  \rangle. \eea

\noindent  Positive definiteness  follows from the fact that,
setting \\$\xi_i(t)=[U_t^{(\underline{\epsilon}_i)}  (
\underbar{u}_i,\underbar{v}_i  ) -\langle
\underbar{u}_i,\underbar{v}_i \rangle  ] \Omega,$ \bean &&
\sum_{i,j=1}^N  \bar{c}_i c_j K \left(
(\underbar{u}_i,\underbar{v}_i,\underline{\epsilon}_i),
 (\underbar{u}_j,\underbar{v}_j,\underline{\epsilon}_j)   \right )\\
&& =\lim_{t\rightarrow 0}\frac{1}{t} \| \sum_{i=1}^N  c_i
\xi_i(t)\|^2\ge 0. \eean

 \noindent {\bf (ii)}  Kolmogorov's
construction (Ref. \cite{krp}) to the pair $(M,K)$  gives the
Hilbert space $\mathbf k$ and embedding $\eta$ satisfying
(\ref{eta-dense}). The separability of $\mathbf k$ follows from
(\ref{eta-dense}) and the verifiable fact $ \| \eta
(\underbar{u},\underbar{v},\underline{\epsilon})-
\eta(\underbar{p},\underbar{w},\underline{\epsilon}) \|_{\mathbf
k}\rightarrow 0$ as
 $  \|\underbar{u}-\underbar{p}\|$ and $\|\underbar{v}-
 \underbar{w}\|\rightarrow 0.$

Setting $\pi$ by (\ref{repn}) we show  that the map $\pi
(\underbar{u},\underbar{v},\underline{\epsilon})$
  extends to a  bounded linear operator on $\mathbf k$ with
$\|\pi (\underbar{u},\underbar{v},\underline{\epsilon})\| \le
\|\underbar{u}\|~\|\underbar{v}\|.$
 For any $\xi=\sum_{i=1}^N  c_i \eta (\underbar{u}_i,\underbar{v}_i,\underline{\epsilon}_i)\in \mathbf k$ let us consider
\bean
&& \|\pi (\underbar{u},\underbar{v},\underline{\epsilon}) \xi\|^2\\
&&=\sum_{i,j=1}^N  \bar{c_i} c_j \langle \pi
(\underbar{u},\underbar{v},\underline{\epsilon}) \eta
(\underbar{u}_i,\underbar{v}_i,\underline{\epsilon}_i),\pi
(\underbar{u},\underbar{v},\underline{\epsilon}) \eta
(\underbar{u}_j,\underbar{v}_j,\underline{\epsilon}_j)
 \rangle\\
&&=\sum_{i,j=1}^N  \bar{c_i} c_j \langle  [\eta (\underbar{u}\otimes \underbar{u}_i,\underbar{v}\otimes\underbar{v}_i,\underline{\epsilon}\oplus \underline{\epsilon}_i)- \langle \underbar{u}_i^{(\underline{\epsilon}_i)},\underbar{v}_i^{(\underline{\epsilon}_i)} \rangle \eta (\underbar{u},\underbar{v},\underline{\epsilon})],\\
&&~~~~~~~~~~~~~~~~~~~[\eta (\underbar{u}\otimes
\underbar{u}_j,\underbar{v}\otimes\underbar{v}_j,\underline{\epsilon}\oplus
\underline{\epsilon}_j)- \langle
\underbar{u}_j^{(\underline{\epsilon}_j)},\underbar{v}_j^{(\underline{\epsilon}_j)}
\rangle \eta (\underbar{u},\underbar{v},\underline{\epsilon})]
 \rangle\\
&&=\lim_{t \rightarrow 0} \frac{1}{t} \sum_{i,j=1}^N  \bar{c_i} c_j \langle U_t^{(\underline{\epsilon})} (\underbar{u}^{(\underline{\epsilon})},\underbar{v}^{(\underline{\epsilon})}) [U_t^{(\underline{\epsilon}_i)}-1] (\underbar{u}_i,\underbar{v}_i) \Omega, U_t^{(\underline{\epsilon})} (\underbar{u}^{(\underline{\epsilon})},\underbar{v}^{(\underline{\epsilon})}) [U_t^{(\underline{\epsilon}_j)}-1] (\underbar{u}_j^{(\underline{\epsilon}_j)},\underbar{v}_j^{(\underline{\epsilon}_j)}) \Omega \rangle\\
&&=\lim_{t \rightarrow 0} \frac{1}{t}  \langle \phi(t), [U_t^{(\underline{\epsilon})} (\underbar{u}^{(\underline{\epsilon})},\underbar{v}^{(\underline{\epsilon})})]^*  U_t^{(\underline{\epsilon})} (\underbar{u}^{(\underline{\epsilon})},\underbar{v}^{(\underline{\epsilon})}) \phi(t) \rangle\\
&&=\lim_{t \rightarrow 0} \frac{1}{t} \|
U_t^{(\underline{\epsilon})}
(\underbar{u}^{(\underline{\epsilon})},\underbar{v}^{(\underline{\epsilon})})
\phi(t) \|^2, \eean where $\phi(t):= \sum_{i,=1}^N   c_i
[U_t^{(\underline{\epsilon}_i)}-1]
(\underbar{u}_i^{(\underline{\epsilon}_i)},\underbar{v}_i^{(\underline{\epsilon}_i)})
\Omega \in \mathcal H.$ In the above identities we have used the
fact that for any $\underline{\epsilon}\in \mathbb Z_2^m,
\underline{\alpha} \in \mathbb Z_2^n$  and product vectors
$\underbar{p}^{(\underline{\epsilon})},\underbar{w}^{(\underline{\epsilon})}
\in \mathbf h^{(\underline{\epsilon})},
\underbar{x}^{(\underline{\alpha})},~\underbar{y}^{(\underline{\alpha})}\in
\mathbf h^{(\underline{\alpha})}$

\bea \label{piUt}
&&[U_t^{\underline{\epsilon}\oplus \underline{\alpha}}-1](\underbar{p}^{(\underline{\epsilon})}\otimes \underbar{x}^{(\underline{\alpha})},\underbar{w}^{(\underline{\epsilon})}\otimes \underbar{y}^{(\underline{\alpha})})- \langle \underbar{x}^{(\underline{\alpha})},\underbar{y}^{(\underline{\alpha})} \rangle [U_t^{(\underline{\epsilon})}-1] (\underbar{p}^{(\underline{\epsilon})},\underbar{w}^{(\underline{\epsilon})}) {\nonumber}\\
&&=U_t^{(\underline{\epsilon})}
(\underbar{p}^{(\underline{\epsilon})},\underbar{w}^{(\underline{\epsilon})})
[U_t^{(\underline{\alpha})}-1]
(\underbar{x}^{(\underline{\alpha})},\underbar{y}^{(\underline{\alpha})}).
\eea Since $ U_t^{(\underline{\epsilon})}
(\underbar{u}^{(\underline{\epsilon})},\underbar{v}^{(\underline{\epsilon})})$
has its norm  bounded  by $\|\underbar{u}\|^2~\| \underbar{v}\|^2$
we get \bean
&& \|\pi (\underbar{u},\underbar{v},\underline{\epsilon}) \xi\|^2\le
\|\underbar{u}\|^2~\| \underbar{v}\|^2 \lim_{t \rightarrow 0} \frac{1}{t}  \|\phi(t)\|^2\\
&&=\sum_{i,j=1}^N   \bar{c_i} c_j \lim_{t \rightarrow 0} \frac{1}{t} \langle   [U_t^{(\underline{\epsilon}_i)}-1] (\underbar{u}_i^{(\underline{\epsilon}_i)},\underbar{v}_i^{(\underline{\epsilon}_i)}) \Omega,[U_t^{(\underline{\epsilon}_j)}-1] (\underbar{u}_j^{(\underline{\epsilon}_j)},\underbar{v}_j^{(\underline{\epsilon}_j)}) \Omega \rangle\\
&& = \|\underbar{u}\|^2~\| \underbar{v}\|^2\|\xi\|^2\eean which
proves that  $\pi (\underbar{u},\underbar{v},\underline{\epsilon})$
extends to a  bounded operator on $\mathbf k$ with\\
$\|\pi (\underbar{u},\underbar{v},\underline{\epsilon})\| \le \|\underbar{u}\|~\|\underbar{v}\|.$\\

 In order to  prove that  $\pi$ is  a $*$-representation of the
algebra $M$  it is enough to  show that for any
$\underline{\epsilon}\in \mathbb Z_2^m, \underline{\epsilon}^\prime
\in \mathbb Z_2^n,\underline{\epsilon}^{\prime \prime} \in \mathbb
Z_2^q$  and product vectors $\underbar{p},\underbar{w} \in \mathbf
h^{\otimes m}, \underbar{p}^\prime,\underbar{w}^\prime \in \mathbf
h^{\otimes n},\underbar{x},\underbar{y} \in \mathbf h^{\otimes q}$

\begin{description}
\item[(i)]$
\pi
(\underbar{u},\underbar{v},\underline{\epsilon})\pi(\underbar{p},\underbar{w},\underline{\epsilon}^\prime)\eta
(\underbar{x},\underbar{y},\underline{\epsilon}^{\prime \prime})=\pi
(\underbar{u}\otimes \underbar{p},\underbar{v}\otimes
\underbar{w},\underline{\epsilon} \oplus
\underline{\epsilon}^\prime)\eta
(\underbar{x},\underbar{y},\underline{\epsilon}^{\prime \prime})$
\item[(ii)]$
\langle \pi
(\underbar{u},\underbar{v},\underline{\epsilon})\eta(\underbar{p},\underbar{w},\underline{\epsilon}^\prime),~\eta
(\underbar{x},\underbar{y},\underline{\epsilon}^{\prime
\prime})\rangle=\langle
\eta(\underbar{p},\underbar{w},\underline{\epsilon}^\prime),~\pi
(\underleftarrow{\mbox{u}},\underleftarrow{\mbox{v}},\underline{\epsilon}^*)\eta
(\underbar{x},\underbar{y},\underline{\epsilon}^{\prime
\prime})\rangle.$
\end{description}
By the definition of $\pi$ \bean
&&\pi (\underbar{u},\underbar{v},\underline{\epsilon})\pi(\underbar{p},\underbar{w},\underline{\epsilon}^\prime)\eta (\underbar{x},\underbar{y},\underline{\epsilon}^{\prime \prime})\\
&&=\pi (\underbar{u},\underbar{v},\underline{\epsilon})[\eta (\underbar{p}\otimes \underbar{x},\underbar{w} \otimes \underbar{y},\underline{\epsilon}^\prime \oplus \underline{\epsilon}^{\prime \prime})-\langle \underbar{x},\underbar{y}\rangle \eta(\underbar{p},\underbar{w},\underline{\epsilon}^\prime)] \\
&&=\eta (\underbar{u}\otimes\underbar{p}\otimes
\underbar{x},\underbar{v}\otimes\underbar{w} \otimes
\underbar{y},\underline{\epsilon} \oplus \underline{\epsilon}^\prime
\oplus \underline{\epsilon}^{\prime \prime})- \langle
\underbar{p}\otimes \underbar{x},\underbar{w} \otimes
\underbar{y}\rangle
\eta(\underbar{u},\underbar{v},\underline{\epsilon})\\
&&-\langle \underbar{x},\underbar{y}\rangle[ \eta(\underbar{u}\otimes\underbar{p},\underbar{v}\otimes\underbar{w},\underline{\epsilon} \oplus\underline{\epsilon}^\prime)-\langle \underbar{p},\underbar{w}\rangle\eta(\underbar{u},\underbar{v},\underline{\epsilon})] \\
&&=\eta (\underbar{u}\otimes\underbar{p}\otimes
\underbar{x},\underbar{v}\otimes\underbar{w} \otimes
\underbar{y},\underline{\epsilon} \oplus \underline{\epsilon}^\prime
\oplus \underline{\epsilon}^{\prime \prime})-\langle
\underbar{x},\underbar{y}\rangle
\eta(\underbar{u}\otimes\underbar{p},\underbar{v}\otimes\underbar{w},\underline{\epsilon}
\oplus\underline{\epsilon}^\prime) \eean and (i) follows. To see
(ii) let us look at  the left hand side. By (\ref{piUt}) \bean
&&\langle \pi (\underbar{u},\underbar{v},\underline{\epsilon})\eta(\underbar{p},\underbar{w},\underline{\epsilon}^\prime),~\eta (\underbar{x},\underbar{y},\underline{\epsilon}^{\prime \prime})\rangle\\
&&=\lim_{t\rightarrow 0}\frac{1}{t} \langle
U_t^{(\underline{\epsilon})}(\underbar{u},\underbar{v})
(U_t^{\underline{\epsilon}^\prime}-1)(\underbar{p},\underbar{w})~\Omega,
(U_t^{\underline{\epsilon}^{\prime \prime}}-1)(\underbar{x},\underbar{y})~\Omega \rangle\\
&&=\lim_{t\rightarrow 0}\frac{1}{t} \langle
(U_t^{\underline{\epsilon}^\prime}-1)(\underbar{p},\underbar{w})~\Omega,
U_t^{\underline{\epsilon}^*}(\underleftarrow{\mbox{v}},\underleftarrow{\mbox{u}})
(U_t^{\underline{\epsilon}^{\prime \prime}}-1)(\underbar{x},\underbar{y})~\Omega \rangle\\
&& =\langle
\eta(\underbar{p},\underbar{w},\underline{\epsilon}^\prime),~\pi
(\underleftarrow{\mbox{u}},\underleftarrow{\mbox{v}},\underline{\epsilon}^*)\eta
(\underbar{x},\underbar{y},\underline{\epsilon}^{\prime
\prime})\rangle=RHS. \eean

Thus \bea \label{pi-mul-*}
&&\pi (\underbar{u},\underbar{v},\underline{\epsilon})\pi(\underbar{p},\underbar{w},\underline{\epsilon}^\prime)=\pi (\underbar{u}\otimes \underbar{p},\underbar{v}\otimes \underbar{w},\underline{\epsilon} \oplus \underline{\epsilon}^\prime) {\nonumber}\\
&&
 \pi (\underbar{u},\underbar{v},\underline{\epsilon})^*=\pi (\underleftarrow{\mbox{u}},\underleftarrow{\mbox{v}},\underline{\epsilon}^*).
\eea
\end{proof}

 \blema
 \label{etauvtotal}
\begin{description}
\item[(a)] For any
$(\underbar{u},\underbar{v},\underline{\epsilon})\in M_0,$
$\underbar{u}=\otimes_{i=1}^n u_i,\underbar{v}=\otimes_{i=1}^n
v_i$\\
and $ \underline{\epsilon}=(\epsilon_1,\epsilon_2,\cdots,
\epsilon_n)$
 \be \label{eta-n}
\eta(\underbar{u},\underbar{v},\underline{\epsilon})= \sum_{ i=1}^n
 \prod_{k\ne i} \langle u_k,v_k \rangle \eta(u_i,v_i, \epsilon_i)
 \ee
\item[(b)] $\eta(u,v,1)=-\eta(u,v,0),~\forall u,v\in \mathbf h.$
\end{description}
 \elema
\begin{proof}

{\bf (a)} For any
$(\underbar{p},\underbar{w},\underline{\epsilon}^\prime)\in M_0,$
 by (\ref{Ut-1-decompose})  and Lemma \ref{4Ut-1},~ we have
\bean && \langle
\eta(\underbar{u},\underbar{v},\underline{\epsilon}),
\eta(\underbar{p}, \underbar{w},\underline{\epsilon}^\prime)
\rangle=K\left( (\underbar{u},\underbar{v},\underline{\epsilon}),
(\underbar{p}, \underbar{w},\underline{\epsilon}^\prime)
\right)\\
&&=\lim_{t\rightarrow 0}\frac{1}{t} \langle
(U_t^{(\underline{\epsilon})}-1)(\underbar{u}, \underbar{v})
\Omega,~ (U_t^{\underline{\epsilon}^\prime}-1)(\underbar{p},
\underbar{w})~\Omega \rangle\\
&&= \sum_{i=1}^n
 \prod_{k\ne i} \overline{\langle u_k,v_k \rangle}
 \lim_{t\rightarrow 0}\frac{1}{t}
\langle (U_t-1)^{(\epsilon_i)}(u_i,v_i)
~\Omega,(U_t^{\underline{\epsilon}^\prime}-1)(\underbar{p},
\underbar{w}) ~\Omega\rangle\\
&&= \sum_{i=1}^n
 \prod_{k\ne i} \overline{\langle u_k,v_k \rangle}
\langle \eta(u_i,v_i,\epsilon_i),\eta(\underbar{p},
 \underbar{w},\underline{\epsilon}^\prime)\rangle.
 \eean
Since $\{\eta(\underbar{p},
\underbar{w},\underline{\epsilon}^\prime):(\underbar{p},
\underbar{w},\underline{\epsilon}^\prime)\in M_0\}$ is a total
subset of $\mathbf k,$   (\ref{eta-n}) follows.\\

\noindent {\bf (b)}
 For any $u,v \in \mathbf h,~(\underbar{p},\underbar{w}
,\underline{\epsilon})\in M_0,$  we have
\bean &&\langle \eta(u,v,1),
\eta(\underbar{p},\underbar{w},\underline{\epsilon}\rangle+\langle
\eta(u,v,0),
\eta(\underbar{p},\underbar{w},\underline{\epsilon}\rangle\\
&&=\lim_{t\rightarrow 0}\frac{1}{t}\langle (U_t+U_t^*-2)(u,v)
\Omega,
(U_t^{(\underline{\epsilon})}-1)(\underbar{p},\underbar{w})~\Omega\rangle\\
&&=-\lim_{t\rightarrow 0}\frac{1}{t}\langle[ (U_t^*-1)(U_t-1)](u,v)
\Omega,
(U_t^{(\underline{\epsilon})}-1)(\underbar{p},\underbar{w})~\Omega\rangle
~~
{\mbox{(since $U_t$ is unitary)}}\\
&&=-\lim_{t\rightarrow 0}\frac{1}{t}
 \sum_{m\ge 1} \langle   (U_t-1)(e_m,u) \Omega,
  (U_t-1)(e_m,v)(U_t^{(\underline{\epsilon})}-1)
  (\underbar{p},\underbar{w})~\Omega\rangle.
\eean

\noindent That this limit  vanishes   can be seen from the following
 \bean &&|
\frac{1}{t}
 \sum_{m\ge 1} \langle   (U_t-1)(e_m,u) \Omega,
  (U_t-1)(e_m,v)(U_t^{(\underline{\epsilon})}-1)
  (\underbar{p},\underbar{w})~\Omega\rangle|^2\\
  &&\le
 \sum_{m\ge 1} \frac{1}{t} \|(U_t-1)(e_m,u) \Omega\|^2
~~ \sum_{m\ge 1} \frac{1}{t}
\|(U_t-1)(e_m,v)(U_t^{(\underline{\epsilon})}-1)
  (\underbar{p},\underbar{w})~\Omega\|^2.\eean
By Lemma  \ref{UTkn} (v) and Lemma  \ref{Auv} (iv) the above
quantity is equal to \bean && 2 Re \langle u, \frac{1-T_t}{t}
u\rangle  \frac{1}{t} \langle   (U_t^{(\underline{\epsilon})}-1)
  (\underbar{p},\underbar{w})~\Omega,
 [(U_t^*-1)(U_t-1)](v,v) (U_t^{(\underline{\epsilon})}-1)
  (\underbar{p},\underbar{w})~\Omega  \rangle\\
&&\le 2 Re \langle u, \frac{1-T_t}{t} u\rangle  \frac{1}{t} \langle
(U_t^{(\underline{\epsilon})}-1)
  (\underbar{p},\underbar{w})~\Omega,
 (2-U_t^*-U_t)(v,v) (U_t^{(\underline{\epsilon})}-1)
  (\underbar{p},\underbar{w})~\Omega  \rangle
\eean Therefore, since $T_t$ is   continuous, by Assumption {\bf C}
\[\lim_{t\rightarrow 0}\frac{1}{t}
 \sum_{m\ge 1} \langle   (U_t-1)(e_m,u) \Omega,
  (U_t-1)(e_m,v)(U_t^{(\underline{\epsilon})}-1)
  (\underbar{p},\underbar{w})~\Omega\rangle=0.\]
 Thus $ \langle \eta(u,v,1), \eta(\underbar{p},\underbar{w},\underline{\epsilon}\rangle
 =-\langle \eta(u,v,0),
 \eta(\underbar{p},\underbar{w},\underline{\epsilon}\rangle.$
As $\{\eta(\underbar{p},\underbar{w},\underline{\epsilon}:
 (\underbar{p},\underbar{w},\underline{\epsilon})\in M_0\}$ is total in $\mathbf k,$
$\eta(u,v,1)=-\eta(u,v,0).$
\end{proof}

\brmrk Writing $\eta(u,v)$ for the vector $\eta(u,v,0)\in \mathbf
k,$ \be \overline{Span}\{ \eta(u,v): u,v\in \mathbf h\}=\mathbf
k.\ee \ermrk

\brmrk The $*$-representation $\pi$ of $M$ in $\mathbf k$ is trivial
\be \pi (\underbar{u},\underbar{v},\underline{\epsilon})
 \eta(\underbar{p},\underbar{w},\underline{\epsilon}^\prime)=
 \langle \underbar{u},\underbar{v}\rangle
 \eta(\underbar{p},\underbar{w},\underline{\epsilon}^\prime) \ee\ermrk
 \noindent
Now we fixed  an ONB    $\{E_j:j\ge 1\}$ for the separable Hilbert
space $\mathbf k.$ Then we have the following  crucial observations.
 \blema \label{H,Lj}
\begin{description}
\item[(a)] There exists a unique family $\{L_j:j\ge 1\}$ in $\mathcal
B(\mathbf h)$ such that $\langle u,L_j v \rangle=\langle E_j,
\eta(u,v) \rangle$ and $
 \sum_{j\ge 1 } \|L_j
u\|^2\le 2 \| G\|\|u\|^2,~ \forall ~u\in \mathbf h,$ so that
$\sum_{j\ge 1 } L_j^*L_j $ converges strongly.
\item[(b)]
The family of operators $\{L_j:j\ge 1\}$  is linearly independent,
 i.e. $\sum_{j\ge 1} c_j L_j=0$  for some $c=(c_j)\in l^2(\mathbb N)$
  implies $c_j=0,~\forall j.$
\item[(c)] If we set  $ iH:=G+\frac{1}{2}
\sum_{j\ge 1} L_j^* L_j$ then $H$  is a bounded self-adjoint
operator on $\mathbf h.$
\end{description}
\elema
\begin{proof}
{\bf(a)} By (\ref{kernel}),   for any $u,v\in \mathbf h$ \bean
&&\|\eta(u,v)\|^2\\
&&= \langle u,\mathcal L (|v><v|) u \rangle-\overline{\langle u,v
\rangle} \langle u,G~v  \rangle- \overline{\langle u, G ~v
\rangle}\langle u,v  \rangle\\
 &&\le  \left[\|\mathcal L\| +2 \| G\|\right
]~\|u\|^2~\|v\|^2. \eean So for each $j\ge 1,$ the map
$\eta_j(u,v):= \langle E_j, \eta(u,v) \rangle,$ defines a bounded
quadratic form on $\mathbf h$ and hence by Riesz's representation
theorem  there exists a unique bounded operator  $L_j \in \mathcal
B(\mathbf h)$ such that $\langle u,L_j v \rangle=\eta_j(u,v).$
 Now  consider the following \bean && \sum_{j} \|L_j
u  \|^2=\sum_{j,k} |\eta_j(e_k,u)|^2
=\sum_{k} \|\eta(e_k,u)\|^2\\
&&=\sum_{k}\left [\langle e_k, \mathcal L(|u><u|)
e_k\rangle-\overline{\langle e_k, u \rangle} \langle e_k, G~ u
\rangle
- \overline{\langle e_k, G~u  \rangle} \langle e_k, u \rangle ~~\right ]\\
&&=Tr \mathcal L(|u><u|)  -\langle u, G~u \rangle-\overline{\langle
u, G~u \rangle}.\eean

\noindent Since $Z_t$ is trace preserving \be \label{lj*lj} \sum_{j}
\|L_j u \|^2 = -\langle u, G~u \rangle-\overline{\langle u, G~u
\rangle}\le 2 \|G\|~\|u\|^2. \ee

\noindent{\bf(b)} Let $\sum_{j\ge 1} c_j L_j=0$  for some
$c=(c_j)\in l^2(\mathbb N).$
 Then for any $u,v\in \mathbf h$ we have
\[0=\langle u, \sum_{j\ge 1} c_j L_j v \rangle = \sum_{j\ge 1} c_j
\langle u, L_j v \rangle =\langle
 \sum_{j\ge 1} \overline{c}_j E_j,\eta (u,v) \rangle.\]
Since $\overline{Span}\{ \eta(u,v): u,v\in \mathbf h\}=\mathbf k ,$
it follows that $
 \sum_{j\ge 1} \overline{c}_j E_j=0\in \mathbf k$  and hence $c_j=0,~\forall j.$

\noindent {\bf(c)} The boundedness of  $G$ and (\ref{lj*lj})
 imply that
$\sum_{j\ge 1} L_j^* L_j$ is a bounded  self-adjoint operator and
hence $H$ is bounded.
  For any   $u\in \mathbf h$  by the identity  (\ref{lj*lj})
\bean
&&\langle  u,( 2G+ \sum_{j\ge 1} L_j^* L_j)u \rangle\\
&&=\langle  u, 2 G u \rangle+ \sum_{j} \|L_ju\|^2
=\langle  u,  G u \rangle-\langle  G   u,  u \rangle\\
&&=-\langle ( 2G+ \sum_{j\ge 1} L_j^* L_j) u,u \rangle \eean Thus
$\langle  u,Hu    \rangle=\langle  H u,u    \rangle$ and by applying
the Polarization  principle to the sesqui-linear form $(u,v)\mapsto
\langle  u,Hu    \rangle $ it  proves that $H$ is self-adjoint.
\end{proof}

\blema The generator  $\mathcal L$ of the uniformly continuous
semigroup $Z_t$  on $\B_1(\mathbf h)$ satisfies \be \label{LZ}
\mathcal L \rho = G\rho+\rho G^*+\sum_{j\ge 1} L_j \rho L_j^*,~~
\forall \rho \in \B_1(\mathbf h). \ee
 \elema
\begin{proof}
By  (\ref{kernel}), for any $u,v,p,w\in \mathbf h$  we have
 \bean
&& \langle \eta(u,v), \eta(p,w) \rangle = \sum_{j\ge 1} \overline{
\langle u, L_j v \rangle } \langle p, L_j w \rangle
 \\
 &&= \langle p,\mathcal L (|w><v|) u \rangle-\overline{\langle u,v
\rangle} \langle p,G~w  \rangle- \overline{\langle u, G ~v
\rangle}\langle p,w  \rangle, \eean which gives \bean
&&\langle p, \mathcal L (|w><v|)~u\rangle \\
&&= \langle p, |Gw><v| ~u\rangle + \langle p,|w><Gv| ~u\rangle
+\sum_{j\ge 1} \langle p, |L_j w>< L_j v|~u\rangle\\
&&= \langle p, G|w><v| ~u\rangle + \langle p,|w><v|G^* ~u\rangle
+\sum_{j\ge 1} \langle p, L_j |w><v| L_j^*~u\rangle . \eean Since
all the operators  involved are bounded (\ref{LZ}) follows.

\end{proof}

\subsection{Associated Hudson-Parthasarathy (HP)  Flows}
Recall from previous section that starting from the family of
unitary operators $\{U_{s,t}\}$  with hypothesis  $\bf A, B,C$  we
obtained the noise  Hilbert  space  $\mathbf k$ and  bounded linear
operators  $G, L_j:j\ge 1$  on the initial Hilbert space  $\mathbf
h.$  Now define a family of operator $\{L_\nu^\mu:\mu, \nu\ge 0\}$
in $\B(\mathbf h)$ by
\begin{equation} \label{hpcoefi}
L_\nu^\mu=\left\{ \begin{array} {lll}
 & G=iH-\half\sum_{k\ge 1}L_k^*L_k & \mbox{for}\
(\mu,\nu)=(0,0)\\
&L_j &  \mbox{for}\ (\mu,\nu)=(j,0)
\\
& -L_k^*&
\mbox{for}\ (\mu,\nu)=(0,k)\\
& 0 & \mbox{for}\ (\mu,\nu)=(j,k).
 \end{array}
 \right.
  \end{equation}
  Note that the indices $\mu, \nu$ vary over  non negative integers while
$j,k$  vary over  non zero positive integers.

\noindent Let us consider the HP type quantum stochastic
differential equation in $\mathbf h \otimes \Gamma(L^2(\mathbb R_+,
\mathbf k))$: \be \label{hpeqn,st} V_{s,t}=1_{\mathbf h \otimes
\Gamma}+\sum_{\mu,\nu\ge 0}\int_s^t V_{s,r} L_\nu^\mu
\Lambda_\mu^\nu(dr)\ee with bounded operator coefficients
$L_\nu^\mu$ given by (\ref{hpcoefi}). By Theorem \ref{hpflow}, there
exists unique  unitary solution $\{V_{s,t}\}$ of the above HP
equation. We shall write  $V_t:=V_{0,t}$ for simplicity. The family
$\{V_{s,t}^*\}$ satisfies: \be \label{hpeqn*}
dV_{s,t}^*=\sum_{\mu,\nu\ge 0} (L_\mu^\nu)^*V_{s,t}^*
\Lambda_\mu^\nu(dt),~ V_{s,s}= 1_{\mathbf h \otimes  \Gamma}\ee and
for any  $u,v\in \mathbf h, V_{s,t}(u,v)$ and $V_{s,t}(u,v)^*$
 satisfy the following  qsde on $\Gamma~:$
\be \label{hpeqn,uv} dV_{s,t}(u,v)=\sum_{\mu,\nu\ge 0}
 V_{s,t}(u,L_\nu^\mu v) \Lambda_\mu^\nu(dt),~ V_{s,s}(u,v)
=\langle u,v \rangle 1_{ \Gamma}.\ee \be \label{hpeqn,*uv}
dV_{s,t}^*(u,v)=\sum_{\mu,\nu\ge 0}
  V_{s,t}^*(L_\mu^\nu u,v) \Lambda_\mu^\nu(dt),~ V_{s,s}^*(u,v)]
=\langle  u, v \rangle1_{ \Gamma}.\ee

\noindent As for the family of unitary operators $\{ U_{s,t}\}$  on
$\mathbf h \otimes \mathcal H,$ for
$\underline{\epsilon}=(\epsilon_1,\epsilon_2,\cdots, \epsilon_n)\in
\mathbb Z_2^n$ we define $V_{s,t}^{(\underline{\epsilon})}\in
\mathcal B(\mathbf h^{\otimes n}\otimes \Gamma)$  by setting
$V_{s,t}^{(\epsilon)}\in \mathcal B(\mathbf h\otimes \Gamma)$ by
\bean
&& V_{s,t}^{(\epsilon)}=V_{s,t} ~\mbox{for}~ \epsilon =0\\
&&~~~~=V_{s,t}^* ~\mbox{for}~ \epsilon =1. \eean We shall write
$V_{s,t}^{(n)}$ for $V_{s,t}^{(\underline{0})},~\underline{0}\in
\mathbb Z_2^n.$

\blema \label{Vstbasic1} The family of unitary operators  $\{
V_{s,t}\}$ satisfy
\begin{description}
\item[(i)] For any $0\le r\le  s \le t<\infty, V_{r,t}=V_{r,s} V_{s,t}.$
\item[(ii)] For $[q,r)\cap [s,t)=\varnothing, V_{q,r}(u,v)$ commute with
 $V_{s,t}(p,w)$  and $V_{s,t}(p,w)^*$  for every $u,v,p,w \in \mathbf h.$
\item[(iii)] For any $0\le s \le t<\infty,$\\
$ \langle \textbf{e}(0),   V_{s,t}(u,v) \textbf{e}(0) \rangle =
\langle \textbf{e}(0),   V_{t-s}(u,v) \textbf{e}(0) \rangle =
\langle u,  T_{t-s }v\rangle,~\forall u,v \in \mathbf h.$
\end{description}
\elema
\begin{proof}
{\bf{(i)}} For fixed $0\le r\le  s\le t<\infty,$  we set
$W_{r,t}=V_{r,s} V_{s,t}$ and $W_{r,s}=V_{r,s} .$
 Then by (\ref{hpeqn,st}) we have
\bean
&& W_{r,t}=V_{r,s}+\sum_{\mu,\nu\ge 0}\int_s^t V_{r,s}V_{s,q} L_\nu^\mu  \Lambda_\mu^\nu(dq)\\
&&=W_{r,s}+\sum_{\mu,\nu\ge 0}\int_s^t   W_{r,q} L_\nu^\mu
\Lambda_\mu^\nu(dq). \eean
 Thus  the family of  unitary operators $\{ W_{r,t}\}$ also satisfies the
 HP equation
(\ref{hpeqn,st}) and, hence by uniqueness of the solution of this
qsde, $W_{r,t}=V_{r,t},
\forall t\ge s$ and the result follows.\\

 \noindent {\bf{(ii)}} For any $0\le s\le t<\infty$
$V_{s,t}\in \mathcal B( \mathbf h\otimes \Gamma_{[s,t]}).$
 So for $p,w\in \mathbf h,~V_{s,t}(p,w)\in \mathcal B(\Gamma_{[s,t]})$ and
the statement  follows.\\

  \noindent {\bf{(iii)}} Let us set  a family of contraction operators $\{
  \widetilde{S}_{s,t}\}$ on $\mathbf h$ by
\[ \langle   u, \widetilde{S}_{s,t} v\rangle =\langle u\otimes \textbf{e}(0),
  V_{s,t}v \otimes \textbf{e}(0) \rangle,~\forall u,v \in \mathbf h.\]
   Then for  fixed $s\ge 0,$  this one parameter family $\{\widetilde{S}_{s,t}\}$
   satisfies the following differential equation
\[ \frac{d \widetilde{S}_{s,t}} {dt}=\widetilde{S}_{s,t} G\] where $G~~
(=L_0^0)$ is the generator of the  uniformly continuous semigroup
$\{T_t\}$ so $\widetilde{S}_{s,t}=T_{t-s}$ and  this proves the
claim.

\end{proof}

 \noindent Consider the  family of maps $\widetilde{Z}_{s,t}$ defined by
\[\widetilde{Z}_{s,t} \rho = Tr_{\mathcal H}  [V_{s,t} (\rho \otimes
| \textbf{e}(0)>< \textbf{e}(0)|) V_{s,t}^*], ~ \forall \rho \in
\B_1(\mathbf h).\]    As for  $Z_t,$ it can be  easily seen
   that $\widetilde{Z}_{s,t}$ is a
  contractive family of maps  on $\B_1(\mathbf h)$  and
  in particular, for any $ u,v,p,w \in \mathbf h$
\[\langle   p, \widetilde{Z}_{s,t} (|w><v|) ~u \rangle
 =\langle   V_{s,t} (u,v) \textbf{e}(0),
  V_{s,t} (p,w) \textbf{e}(0) \rangle.\]
  \blema
  The family $\widetilde{Z}_t:=\widetilde{Z}_{0,t} $  is a  uniformly continuous
  semigroup  of contraction  on $\B_1(\mathbf h)$  and
  $\widetilde{Z}_{s,t}=\widetilde{Z}_{t-s}
  =Z_{t-s}.$
  \elema
  \begin{proof}  By (\ref{hpeqn,uv}) and Ito's formula
\bean &&\langle   p,[ \widetilde{Z}_{s,t}-1] (|w><v|) ~u \rangle\\
  &&=\langle   V_{s,t} (u,v) \textbf{e}(0),
  V_{s,t} (p,w) \textbf{e}(0) \rangle-\overline{\langle   u,v
   \rangle} \langle  p,w \rangle\\
   &&=\int_s^t \langle   V_{s,\tau} (u,v) \textbf{e}(0),
  V_{s,\tau} (p,G w) \textbf{e}(0) \rangle d\tau
  +\int_s^t \langle   V_{s,\tau} (u,Gv) \textbf{e}(0),
  V_{s,\tau} (p,w) \textbf{e}(0) \rangle d\tau\\
 && +\int_s^t
  \langle   V_{s,\tau} (u,L_jv) \textbf{e}(0),
  V_{s,\tau} (p,L_j w) \textbf{e}(0) \rangle d\tau\\
  &&=\int_s^t \langle   p, \widetilde{Z}_{s,\tau} (|Gw>< v|) ~u \rangle
  d\tau
  +\int_s^t\langle   p, \widetilde{Z}_{s,\tau} (|w><G v|) ~u \rangle  d\tau\\
 &&+\sum_{j\ge 1}\int_s^t\langle   p, \widetilde{Z}_{s,\tau}
  (|L_jw><L_j v|) ~u \rangle
 d\tau\\
 && =\int_s^t \langle   p, \widetilde{Z}_{s,\tau}
 \mathcal L(|w>< v|) ~u \rangle d\tau,
  \eean where $\mathcal L $  is the generator of the  uniformly continuous
  semigroup $Z_t.$
 Since the maps  $\mathcal L$ and  $\widetilde{Z}_{a,b}:0\le a\le b$ are  bounded,
 for fixed $s\ge 0,~ \widetilde{Z}_{s,t}$ satisfies the differential equation
\[\widetilde{Z}_{s,t} (\rho) =\rho+\int_s^t \widetilde{Z}_{s,\tau}
 \mathcal L(\rho)  d\tau,~~\rho\in  \B_1(\mathbf h) .\] Hence $\widetilde{Z}_t $  is a
  uniformly continuous
  semigroup   on $\B_1(\mathbf h)$  and
$\widetilde{Z}_{s,t}=\widetilde{Z}_{t-s}
  =Z_{t-s}.$
  \end{proof}

\section{Minimality of HP Flows}
\noindent In this section  we shall show the minimality of the HP
flow $V_{s,t}$ discussed above. We prove that the subset $ \mathcal
S^\prime:=\{ \zeta=
V_{\underbar{s},\underbar{t}}(\underbar{u},\underbar{v})
\textbf{e}(0):=V_{s_1,t_1}(u_1,v_1)\cdots
V_{s_n,t_n}(u_n,v_n)\textbf{e}(0): \underbar{s}=(s_1,s_2, \cdots,
s_n), \underbar{t}=(t_1,t_2, \cdots, t_n)$ $:~ 0 \le s_1\le t_1\le
s_2 \le \ldots\le  s_n\le  t_n< \infty,n\ge 1,
\underbar{u}=\otimes_{i=1}^n u_i,\underbar{v}=\otimes_{i=1}^n v_i\in
\mathbf h^{\otimes n}\}$  is total in  the symmetric Fock space
$\Gamma(L^2(\mathbb R_+,\mathbf k)).$

\noindent  We note that for any $0\le s < t\le \tau<\infty , u,v\in
\mathbf h$ by the HP equation (\ref{hpeqn,st}) \bea \label{vst-1}
&&\frac{1}{t-s}[V_{s,t}-1] (u,v) \textbf{e}(0) {\nonumber}\\
&& = \frac{1}{t-s}\{\sum_{j\ge 1} \int_s^t V_{s,\lambda} (u,L_jv)
a_j^\dag (d\lambda)+ \int_s^t V_{s,\lambda} (u,Gv)
d\lambda  \}\textbf{e}(0) {\nonumber}\\
&& = \gamma(s,t,u,v)+  \langle u,Gv\rangle
~~\textbf{e}(0)+\zeta(s,t,u,v) +\varsigma(s,t,u,v), \eea \noindent
here these vectors in the Fock space $\Gamma$ are given by
\begin{description}
\item $\gamma(s,t,u,v):=\frac{1}{t-s}\sum_{j\ge 1} \langle u,L_jv\rangle
a_j^\dag([s,t])~~\textbf{e}(0)$
\item $\zeta(s,t,u,v):=\frac{1}{t-s}\sum_{j\ge 1}\int_s^t (V_{s,\lambda}
-1)(u,L_jv) a_j^\dag (d\lambda) ~~\textbf{e}(0)$
\item  $\varsigma(s,t,u,v):=\frac{1}{t-s} \int_s^t (V_{s,\lambda}-1)
(u,Gv) d \lambda ~~\textbf{e}(0).$
\end{description}

\noindent Note that any $\xi\in \Gamma $ can be written as
$\xi=\xi^{(0)} \textbf{e}(0)\oplus \xi^{(1)} \oplus \cdots,~
\xi^{(n)}$ is in the $n$-fold symmetric tensor product  $L^2(\mathbb
R_+,\mathbf k)^{\otimes n}\equiv L^2(\Sigma_n) \otimes \mathbf
k^{\otimes n}),$ where $\Sigma_n$ is the $n$-simplex $\{
\underbar{t}=(t_1,t_2,\cdots, t_n):0\le t_1 <t_2< \ldots <t_n
<\infty\}.$

\blema \label{qs-esti} Let $\tau \ge 0.$ For any $ u,v\in \mathbf h,
0\le s\le t\le \tau,$  define constants $C_\tau=2 e^{\tau}$  and $
C_{\tau,v}=  C_\tau \{ \sum_{j\ge 1}
  \|L_j v\|^2 +  \tau\|G~v\|^2\}.$ Then\\

{\bf a.} \be \|(V_{s,t}-1) v \textbf{e}(0)\|^2\le C_{\tau,v}
(t-s).\ee

{\bf b.} For any $u\in \mathbf h$ \bean &&\|\sum_{j\ge 1} \int_s^t
V_{s,\lambda} (u,L_jv) a_j^\dag (d\lambda) \textbf{e}(0)\|^2\\
&&\le C_\tau \| u\|^2  \sum_{j\ge 1}  \int_s^t
\|V_{s,\lambda} L_jv \otimes  \textbf{e}(0) \|^2 ~d\lambda\\
&& \le C_\tau (t-s)\| u\|^2 \sum_{j\ge 1} \|L_jv \|^2.\eean \elema
\begin{proof}{\bf a.} By
estimates of   quantum stochastic integration  (Proposition 27.1,
\cite{krp})
\bean
&&\|(V_{s,t}-1) v \textbf{e}(0)\|^2\\
&&=\| \sum_{j\ge 1}\int_s^t V_{s,\lambda} L_j a_j^\dag (d\lambda)
~~v \textbf{e}(0)+ \int_s^t V_{s,\lambda}G d
\lambda ~~v \textbf{e}(0)\|^2\\
&&\le C_\tau \int_s^t\{  \sum_{j\ge 1}\|L_jv\|^2 +\|Gv\|^2\} d
\lambda \\
 && = C_{\tau,v} (t-s).\eean
 {\bf b.} For  any  $\phi$ in the Fock space $\Gamma(L^2(\mathbb R_+, \mathbf k))$,
\bean && \langle \phi,\sum_{j\ge 1} \int_s^t V_{s,\lambda} (u,L_jv)
a_j^\dag (d\lambda) \textbf{e}(0)\rangle |^2\\
&&=|\langle u\otimes\phi,\{ \sum_{j\ge 1} \int_s^t V_{s,\lambda}L_j
a_j^\dag (d\lambda) \}v  \textbf{e}(0)\rangle |^2\\
&&\le \| u\otimes\phi\|^2 \|\{ \sum_{j\ge 1} \int_s^t
V_{s,\lambda}L_j a_j^\dag (d\lambda)\} v  \textbf{e}(0) \|^2. \eean
By estimates of quantum stochastic integration   the above quantity
is \[\le C_\tau \| u\otimes\phi\|^2  \sum_{j\ge 1}  \int_s^t
\|V_{s,\lambda} L_jv  \textbf{e}(0) \|^2 ~d\lambda.\] Since $\phi$
is arbitrary and  $V_{s,\lambda}$'s are contractive the statement
follows.

\end{proof}
 \blema \label{lemma-vst-1} Let $\tau \ge 0.$ For any $ u,v\in
\mathbf h, 0\le s\le t\le \tau$
\begin{description}
\item[(a)] $\|(V_{s,t}-1)(u,v)~~\textbf{e}(0)\|^2\le 2 C_{\tau,v}\|u\|^2
(t-s).$

\item[(b)] $\sup\{ \|\zeta(s,t,u,v) \|^2:0\le s\le t \le \tau\} <\infty$
and\\
 $\|\varsigma(s,t,u,v)\|\le \|u\| \sqrt{ 2C_{\tau,Gv}(t-s)},~ \forall ~0\le s < t\le
\tau.$

\item[(c)]For any $\xi\in \Gamma(L^2(\mathbb R_+, \mathbf k)),$
~ $\lim_{s \rightarrow t}\langle \xi,  \zeta(s,t,u,v)\rangle =0$ and
\[\lim_{s \rightarrow t} \langle \xi, \gamma(s,t,u,v)\rangle =\sum_{j\ge 1}
\langle u,L_jv\rangle \overline{\xi^{(1)}_j} (t)= \langle \xi^{(1)}
(t), \eta(u,v)\rangle,~~\mbox{a.e.}~~ t\ge 0.\]
\end{description}
\elema

\begin{proof} {\bf (a)}  By identity  (\ref{vst-1}) and Lemma \ref{qs-esti} (b) we have
   \bean
  && \|(V_{s,t}-1)(u,v)~~\textbf{e}(0)\|^2 \\
   &&=\|\sum_{j\ge 1} \int_s^t
  V_{s,\alpha}(u,L_j v) a_j^\dag(d\alpha)~\textbf{e}(0)+ \int_s^t
  V_{s,\alpha}(u, G v)~\textbf{e}(0) d\alpha\|^2\\
  && \le 2 \|\sum_{j\ge 1} \int_s^t
  V_{s,\alpha}(u,L_j v) a_j^\dag(d\alpha)~~\textbf{e}(0)\|^2+[\int_s^t
  \|V_{s,\alpha}(u, G v)~~\textbf{e}(0)\| d\alpha]^2\\
   && \le 2 \|u\|^2[ C_\tau (t-s) \sum_{j\ge 1}
  \|L_j v\|^2 +  [(t-s)\|G~v\|]^2 ] \\
     && \le 2 C_{\tau,v}\|u\|^2  (t-s).
  \eean

\noindent {\bf(b) }  {1.}  As in the  Lemma \ref{qs-esti} (b) we
have
  \bean
  && \|\zeta(s,t,u,v)\|^2=\frac{1}{(t-s)^2}\|\sum_{j\ge 1} \int_s^t
  (V_{s,\lambda}-1)(u,L_jv) a_j^\dag (d\lambda) ~~\textbf{e}(0)\|^2\\
    && \le  \frac{ \|u\|^2}{(t-s)^2}\| \sum_{j\ge 1} \int_s^t
 \|(V_{s,\lambda}-1)L_jv ~~\textbf{e}(0)\|^2 d\lambda.
 \eean
   Now by Lemma \ref{qs-esti} (a),  the above quantity is
 \bean
     && \le  \frac{C_\tau \|u\|^2}{(t-s)^2}\sum_{j\ge 1}
      C_\tau (t-s)^2\{ \sum_{i\ge 1}
  \|L_i L_j v\|^2 +  \tau \|G~L_j ~v\|]^2 \} \\
  && \le  C_\tau^2 \|u\|^2 \{\sum_{j\ge 1}
        \sum_{i\ge 1}
  \|L_i L_j v\|^2 + \tau \sum_{j\ge 1} \|G~L_j ~v\|^2 \}. \eean
     Since $\sum_{j\ge 1} \|L_j ~v\|^2=-2 Re \langle v,Gv \rangle,$  the  above  quantity  is  bounded
     and is independent of $s,t.$

\noindent {2.} We have
  \bean
  && \|\varsigma(s,t,u,v)\|=\frac{1}{(t-s)}\| \int_s^t
  (V_{s,\lambda}-1)
(u,Gv) d \lambda ~~\textbf{e}(0)\|\\
    && \le  \frac{1}{(t-s)} \int_s^t
  \|(V_{s,\lambda}-1)
(u,Gv)~~\textbf{e}(0)\| d\lambda.
  \eean
By {\bf(a)} the estimate follows.

\noindent  {\bf(c)}
 {1.} For any
$f\in L^2(\mathbb R_+,\mathbf k).$
 Let us consider
 \bean
 &&\langle
\textbf{e}(f),\zeta(s,t,u,v)\rangle=\langle
\textbf{e}(f),\frac{1}{t-s} \sum_{j\ge 1}\int_s^t (V_{s,\lambda}-1)
(u,L_jv)
a_j^\dag (d\lambda) ~~ \textbf{e}(0)\rangle\\
&&=\frac{1}{t-s} \sum_{j\ge 1}\int_s^t \overline{f_j(\lambda)}
\langle \textbf{e}(f),(V_{s,\lambda}-1) (u,L_jv)  ~~
\textbf{e}(0)\rangle  d\lambda\\
&&=\frac{1}{t-s} \int_s^t G(s,\lambda) d\lambda, \eean where
$G(s,\lambda)= \sum_{j\ge 1} \overline{f_j(\lambda)} \langle
\textbf{e}(f),(V_{s,\lambda}-1) (u,L_jv)  ~~ \textbf{e}(0)\rangle. $
 Note that the complex valued function
$G(s,\lambda)$ is uniformly continuous in both the variables
$s,\lambda$ on $[0,\tau]$  and $G(t,t)=0.$ So we get
\[\lim_{s \rightarrow t} \langle
\textbf{e}(f),\zeta(s,t,u,v)\rangle=0.\] Since  $\zeta(s,t,u,v)$ is
uniformly bounded in $s,t$
\[\lim_{s \rightarrow t} \langle
\xi,\zeta(s,t,u,v)\rangle=0, \forall \xi \in \Gamma.\]

{2.} We have
 \be \label{xi1}  \langle \xi,
\gamma(s,t,u,v)\rangle = \frac{1}{t-s}\sum_{j\ge 1} \langle
u,L_jv\rangle \int_s^t \overline{\xi^{(1)}_j} (\lambda )
d\lambda.\ee Since
\[|\sum_{j\ge 1} \langle
u,L_jv\rangle \overline{\xi^{(1)}_j} (t) |^2\le \|u\|^2 \sum_{j\ge
1} \|L_jv\|^2 |\xi^{(1)}_j (t)|^2\le C \|v\|^2  \|\xi^{(1)}
(t)\|^2,\] the function $\sum_{j\ge 1} \langle u,L_jv\rangle
\overline{\xi^{(1)}_j} (\cdot)$ is in $ L^2$  and  hence locally
integrable. Thus we get
\[\lim_{s \rightarrow t} \langle \xi,
\gamma(s,t,u,v)\rangle=\sum_{j\ge 1} \langle u,L_jv\rangle
\overline{\xi^{(1)}_j} (t)~~\mbox{a.e.}~~ t\ge 0.\]

\end{proof}
\blema \label{Mst}
 For $n\ge 1, ~ \underbar{t}\in \Sigma_n$  and $u_k,v_k\in
\mathbf h:k=1,2,\cdots, n, \xi \in \Gamma(L^2(\mathbb R_+,\mathbf
k))$ and  disjoint intervals $[s_k,t_k),$
\begin{description}
\item[(a)] $\lim_{\underbar{s} \rightarrow \underbar{t}}\langle \xi,
\prod_{k=1}^n M(s_k,t_k,u_k,v_k)~~\textbf{e}(0) \rangle=0,$ \\
where
$M(s_k,t_k,u_k,v_k)=\frac{(V_{s_k,t_k}-1)}{t_k-s_k}(u_k,v_k)-\langle
u_k, G~v_k \rangle-\gamma(s_k,t_k,u_k,v_k)$  and $\lim_{\underbar{s}
\rightarrow \underbar{t}}$  means  $s_k\rightarrow t_k$  for each
$k.$

\item[(b)]
$\lim_{\underbar{s} \rightarrow \underbar{t}}\langle \xi,
\otimes_{k=1}^n \gamma(s_k,t_k,u_k,v_k) \rangle=\langle
\xi^{(n)}(t_1,t_2,\cdots, t_n), \eta(u_1,v_1) \otimes\cdots \otimes
\eta(u_n,v_n) \rangle.$
\end{description}
 \elema
 \begin{proof}
 {\bf(a)} First note that
$M(s,t,u,v)\textbf{e}(0)=\zeta(s,t,u,v)+~ \varsigma(s,t,u,v).$
 So by the above observations
 $\{M(s,t,u,v) \textbf{e}(0)\}$ is  uniformly bounded in $s,t$
 and\\
$\lim_{s \rightarrow t} \langle \textbf{e}(f)
,M(s,t,u,v)\textbf{e}(0)\rangle=0, \forall f\in L^2(\mathbb
R_+,\mathbf k).$ Since the intervals $[s_k,t_k)$'s are disjoint for
different $k$'s,
\[\langle \textbf{e}(f), \prod_{k=1}^n M(s_k,t_k,u_k,v_k)~~\textbf{e}(0)
\rangle= \prod_{k=1}^n \langle \textbf{e}(f_{[s_k,t_k)}),
M(s_k,t_k,u_k,v_k)~~\textbf{e}(0) \rangle\] and thus
$\lim_{\underbar{s} \rightarrow \underbar{t}}\langle \textbf{e}(f),
\prod_{k=1}^n M(s_k,t_k,u_k,v_k)~~\textbf{e}(0) \rangle=0.$ By Lemma
\ref{lemma-vst-1}, the vector  $\prod_{k=1}^n
M(s_k,t_k,u_k,v_k)~~\textbf{e}(0)$ is uniformly bounded  in
$s_k,t_k$
and the convergence can be extended  to  Fock Space.\\

\noindent {\bf(b) } It can be proved similarly as part {\bf(c) }  of
the  previous Lemma.

\end{proof}

\blema Let $\xi\in \Gamma$ be such that  \be \label{ortho} \langle
\xi, \zeta\rangle=0,~\forall \zeta \in \mathcal S^\prime, \ee Then
\begin{description}
\item[(i)]$\xi^{(0)}=0.$
\item[(ii)] $\xi^{(1)}(t)=0,$
~\mbox{for  a.e.} $t\in[0,\tau].$
\item[(iii)]For any $n\ge 0,~ \xi^{(n)}( \underbar{t})=0,$
~\mbox{for  a.e.} $\underbar{t}\in \Sigma_n: t_i\le \tau.$
\item[(iv)]The set  $\mathcal S^\prime$ is total in the  Fock space $\Gamma.$
\end{description}
\elema

\begin{proof}
{\bf(i)} For any $s\ge 0,~V_{s,s}=1_{\mathbf h\otimes \Gamma}$ so in
particular (\ref{ortho}) gives, for any  $u,v\in \mathbf h$
\[0=\langle \xi, V_{s,s}(u,v)  \textbf{e}(0)\rangle=\langle u,v\rangle
\overline{\xi^{(0)}}\] and hence $\xi^{(0)}=0.$\\

\noindent {\bf(ii)} By (\ref{ortho}),~ $\langle \xi, [V_{s,t}-1]
(u,v) \textbf{e}(0)\rangle=0$ for any  $0\le s < t\le \tau<\infty ,
u,v\in \mathbf h.$ By HP equation (\ref{hpeqn,st}) and
 Lemma \ref{vst-1} we have \bean &&0=\lim_{s \rightarrow t}
\frac{1}{t-s}\langle \xi, [V_{s,t}-1] (u,v) \textbf{e}(0)\rangle \\
&&=\sum_{j\ge 1}\langle u,L_jv \rangle \overline{\xi^{(1)}_j(t)}
=\sum_{j\ge 1}\eta_j(u,v) \overline{\xi_j^{(1)}(t)}. \eean So
$\langle \xi^{(1)}(t), \eta(u,v) \rangle=0, \forall u,v\in \mathbf
h.$ Since $\{\eta(u,v):u,v\in \mathbf h\}$ is total in $\mathbf k$
it follows that $\xi^{(1)}(t)=0$ for  $0\le t\le \tau.$

\noindent {\bf(iii)} We prove this by induction. The result is
already proved for $n=0,1.$ For $n\ge 2,$ assume as induction
hypothesis  that for all $m\le n-1,$ $\xi^{(m)}( \underbar{t})=0,$
~\mbox{for a.e.} $\underbar{t}\in \Sigma_m: t_i\le
\tau,i=1,2,\cdots, m.$ We now   show that $\xi^{(n)}(
\underbar{t})=0,$ ~\mbox{for
a.e.} $\underbar{t}\in \Sigma_n: t_i\le \tau.$\\

 \noindent Let $0\le s_1< t_1\le s_2<t_2<\ldots < s_n <t_n\le \tau$
 and $u_i,v_i\in \mathbf h:i=1,2\cdots, n.$
By (\ref{ortho})  and  part ({\bf i})  we have
\[\langle \xi,
 \prod_{k=1}^n \frac{(V_{s_k,t_k}-1)}{t_k-s_k}(u_k,v_k)~~\textbf{e}(0) \rangle =0.\]
Thus
 \bea \label{nvst-1}
 && 0=\lim_{\underbar{s}
\rightarrow \underbar{t}}\langle \xi,
 \prod_{k=1}^n \frac{(V_{s_k,t_k}-1)}{t_k-s_k}(u_k,v_k)~~\textbf{e}(0) \rangle \\
 && =\lim_{\underbar{s}
\rightarrow \underbar{t}}\langle \xi, \prod_{k=1}^n \{
M(s_k,t_k,u_k,v_k)+\langle u_k, G~v_k \rangle+
 \gamma(s_k,t_k,u_k,v_k) \}~~\textbf{e}(0) \rangle. {\nonumber}
\eea Let $P,Q,R$   and $ P^\prime, R^\prime $ be two sets of
disjoint  partitions of $ \{1,2,\cdots , n\}$ such that $Q$ and $ R$
are non empty. We write $|S|$ for the cardinality  of set  $S.$ Then
by Lemma \ref{Mst} (b) the  right hand side of (\ref{nvst-1}) is
equal to
 \bean
 && \sum_{P^\prime
,R^\prime }  \langle
\xi^{(|R^\prime|)}(t_{r_1^\prime},\cdots,t_{r_{|R^\prime|}^\prime
}), \otimes_{k \in R^\prime }\eta(u_k,v_k) \rangle ~\prod_{k\in
P^\prime
} \langle u_k, G~v_k \rangle\\
 &&~~~ +  \lim_{\underbar{s}
\rightarrow \underbar{t}}\sum_{P,Q,R}\langle \xi,  \prod_{k\in P}
\langle u_k, G~v_k \rangle~\prod_{k\in Q}\{ M
(s_k,t_k,u_k,v_k)\}\prod_{k \in R} \{\gamma(s_k,t_k,u_k,v_k)
\}~~\textbf{e}(0) \rangle.
 \eean
Thus by the induction hypothesis,
 \bea \label{lastterm}&& 0= \langle  \xi^{(n)}(t_1,t_2,\cdots, t_n),
\eta(u_1,v_1)
\otimes\cdots \otimes  \eta(u_n,v_n) \rangle\\
 &&~~~ +  \lim_{\underbar{s}
\rightarrow \underbar{t}}\sum_{P,Q,R}\langle \xi,  \prod_{k\in P}
\langle u_k, G~v_k \rangle~\prod_{k\in Q}\{ M
(s_k,t_k,u_k,v_k)\}\prod_{k \in R} \{\gamma(s_k,t_k,u_k,v_k)
\}~~\textbf{e}(0) \rangle. {\nonumber}
 \eea
We claim that the second  term in  (\ref{lastterm}) vanishes. To
prove the claim, it is enough to  show  that for any two non empty
disjoint subsets $Q\equiv\{q_1,q_2,\cdots,
q_{|Q|}\},R\equiv\{r_1,r_2, \cdots, r_{|R|}\}$  of   $ \{1,2,\cdots,
n\},$
 \be \label{QR00}
 \lim_{\underbar{s}
\rightarrow \underbar{t}} \langle \xi, \prod_{q\in Q}\{ M
(s_q,t_q,u_q,v_q)\}\prod_{r \in R} \{\gamma(s_r,t_r,u_r,v_r)
\}~~\textbf{e}(0) \rangle=0.\ee Writing $ \psi $ for the vector
$\prod_{q\in Q}\{ M (s_q,t_q,u_q,v_q) \}\textbf{e}(0),$ we have
 \bea \label{QR11}
&& \langle \xi, \prod_{q\in Q}\{ M (s_q,t_q,u_q,v_q)\}\prod_{r \in
R} \{\gamma(s_r,t_r,u_r,v_r)
\}~~\textbf{e}(0) \rangle    {\nonumber}\\
&&=  \langle \xi, \psi \otimes  \otimes_{r \in
R} \frac{1_{[s_r,t_r]}\eta(u_r,v_r)}{t_r-s_r} \rangle   {\nonumber}\\
&&=\langle \xi,\psi \otimes \otimes_{r \in R}
\frac{1_{[s_r,t_r]}\eta(u_r,v_r)}{t_r-s_r}  \rangle {\nonumber}\\
&&=  \sum_{l\ge |R|} \langle \xi^{(l)},\psi^{(l-|R|)} \otimes
\otimes_{r \in R} \frac{1_{[s_r,t_r]}\eta(u_r,v_r)}{t_r-s_r}
\rangle {\nonumber}\\
 &&=\langle \sum_{l\ge |R|} \langle \langle
\xi^{(l)},\psi^{(l-|R|)} \rangle \rangle, \otimes_{r \in R}
\frac{1_{[s_r,t_r]}\eta(u_r,v_r)}{t_r-s_r}  \rangle. \eea
 Here  $ \langle \langle
\xi^{(l)},\psi^{(l-|R|)} \rangle \rangle \in L^2(\mathbf R_+,\mathbf
k)^{\otimes |R|}$ is  defined  as in (\ref{partinn}) by
 \bea
\label{<<>>} && \langle ~~\langle \langle \xi^{(l)},\psi^{(l-|R|)}
\rangle \rangle, \rho^{(|R|)}~~ \rangle =\langle
\xi^{(l)},\psi^{(l-|R|)}
\otimes  \rho^{(|R|)} \rangle\\
&&=\int_{\Sigma_l}\langle \xi^{(l)} (x_1,x_2,\cdots,
x_l),\psi^{(l-|R|)} (x_1, x_2,\cdots, x_{l-|R|})\otimes
\rho^{(|R|)} (x_{l-|R|+1},\cdots, x_l) \rangle_{\mathbf k^{\otimes
l}}~~ dx {\nonumber} \eea for any
$\rho^{(|R|)}\in L^2(\mathbf R_+,\mathbf k)^{\otimes |R|}.$\\
 By
Lemma \ref{Mst} (a),
 \be
  \lim_{s_q \rightarrow t_q}
  \langle \xi, \prod_{q\in Q}\{ M
(s_q,t_q,u_q,v_q)\}\prod_{r \in R} \{\gamma(s_r,t_r,u_r,v_r)
\}~~\textbf{e}(0) \rangle=0.\ee  However, we need to prove
(\ref{QR00}) where the limit $\underbar{s} \rightarrow \underbar{t}$
has to be in arbitrary order. On the other hand, by (\ref{QR11}) and
(\ref{<<>>}) we get
 \bea \label{QR33} && \lim_{s_q \rightarrow t_q} \lim_{s_r \rightarrow t_r}\langle \xi,
\prod_{q\in Q}\{ M (s_q,t_q,u_q,v_q)\}\prod_{r \in R}
\{\gamma(s_r,t_r,u_r,v_r) \}~~\textbf{e}(0) \rangle
{\nonumber}\\
&&= \lim_{s_q \rightarrow t_q} \lim_{s_r \rightarrow t_r} \langle
\sum_{l\ge |R|} \langle \langle \xi^{(l)},\psi^{(l-|R|)} \rangle
\rangle, \otimes_{r
\in R} \frac{1_{[s_r,t_r]}\eta(u_r,v_r)}{t_r-s_r}  \rangle{\nonumber}\\
&&= \lim_{s_q \rightarrow t_q} \lim_{s_r \rightarrow t_r} \langle
\int_{\Sigma_{|R|}}\langle [\sum_{l\ge |R|}  \langle \langle
\xi^{(l)},\psi^{(l-|R|)} \rangle \rangle ](x_1,x_2,\cdots, x_{|R|}), {\nonumber}\\
&&~~~~~~~~~~~~~~~~~~~~~~~ \otimes_{r \in R} \frac{1_{[s_r,t_r]}
(x_r)~~\eta(u_r,v_r)}{t_r-s_r} \rangle dx{\nonumber}\\
&&= \lim_{s_q \rightarrow t_q}   \langle \sum_{l\ge |R|} \langle
\langle \xi^{(l)},\psi^{(l-|R|)} \rangle \rangle(t_{r_1}, \cdots,
t_{r_{|R|}}), \otimes_{r \in R}\eta(u_r,v_r) \rangle,\eea for almost
all $\underbar{t}\in \Sigma_{|R|}.$ We fix $\underbar{t}\in
\Sigma_{|R|}$ and define families of vectors  $\widetilde{\xi}^{(l)}
:l\ge 0$  in $ L^2(\mathbb R_+, \mathbf k) ^{\otimes l}$ by
 \bean && \widetilde{\xi}^{(0)} = \langle \xi^{(|R|)}(t_{r_1}, \cdots,
t_{r_{|R|}}), \otimes_{r \in
R}\eta(u_r,v_r) \rangle \in \mathbb C\\
&& \widetilde{\xi}^{(l)} (x_1,x_2,\cdots , x_l) = \langle\langle
\xi^{(|R|+l)}(x_1,\cdots ,x_l, t_{r_1}, \cdots, t_{r_{|R|}}),
\otimes_{r \in R}\eta(u_r,v_r) \rangle \rangle, \eean which defines
a Fock space vector $\widetilde{\xi}.$ Therefore, from (\ref{QR33}),
we get  that \bean && \lim_{s_q \rightarrow t_q} \lim_{s_r
\rightarrow t_r} \langle \xi, \prod_{q\in Q}\{ M
(s_q,t_q,u_q,v_q)\}\prod_{r \in R} \{\gamma(s_r,t_r,u_r,v_r)
\}~~\textbf{e}(0) \rangle = \lim_{s_q
\rightarrow t_q}   \langle \widetilde{\xi}~,~\psi \rangle\\
&&= \lim_{s_q \rightarrow t_q}   \langle \widetilde{\xi}~,~
[\prod_{q\in Q} M (s_q,t_q,u_q,v_q)] ~\textbf{e}(0)\rangle,\eean
which is equal to $0$ by Lemma \ref{Mst}  (a).
 Thus from (\ref{lastterm})  we get that  \[\langle
\xi^{(n)}(t_1,t_2,\cdots, t_n), \eta(u_1,v_1) \otimes\cdots \otimes
\eta(u_n,v_n) \rangle=0.\] Since $\{ \eta(u,v):u,v\in \mathbf h\}$
is total in $\mathbf k,$ it follows that $\xi^{(n)}(t_1,t_2,\cdots,
t_n)=0$ for almost every $(t_1,t_2,\cdots, t_n)\in \Sigma_n:t_k\le
\tau.$
\end{proof}

\noindent {\bf(iv)} Since $\tau\ge 0 $  is arbitrary $\xi^{(n)}=0\in
L^2(\mathbb R_+, \mathbf k)^{\otimes n}:n\ge 0$  and hence $\xi=0.$
Which proves the  totality  of $\mathcal S^\prime\subseteq  \Gamma.$

\section{Unitary Equivalence}
Here  we shall show that  the  unitary evolution $\{U_{s,t}\}$ on
$\mathbf h \otimes \mathcal H $ is unitarily equivalent to the HP
flow $\{V_{s,t}\}$ on $\mathbf h \otimes \Gamma(L^2(\mathbb
R_+,\mathbf k)) $ discussed above.  Let us recall that the subset $
\mathcal S=\{ \xi=
U_{\underbar{s},\underbar{t}}(\underbar{u},\underbar{v})
\Omega:=U_{s_1,t_1}(u_1,v_1)\cdots U_{s_n,t_n}(u_n,v_n)\Omega :
\underbar{s}=(s_1,s_2, \cdots, s_n), \underbar{t}=(t_1,t_2, \cdots,
t_n)$ $:~ 0 \le s_1\le t_1\le  s_2 \le \ldots \le s_n\le  t_n<
\infty,n\ge 1, \underbar{u}=\otimes_{i=1}^n
u_i,\underbar{v}=\otimes_{i=1}^n v_i\in \mathbf h^{\otimes n}\}$ is
total in $\mathcal H$
 and the subset\\
$ \mathcal S^\prime:=\{ \zeta=
V_{\underbar{s},\underbar{t}}(\underbar{u},\underbar{v})
\textbf{e}(0):=V_{s_1,t_1}(u_1,v_1)\cdots
V_{s_n,t_n}(u_n,v_n)\textbf{e}(0):\\~~~~~~~~~~~~~~~~~~~~~
\underbar{u},\underbar{v}\in \mathbf h^{\otimes
n},\underbar{s}=(s_1,s_2, \cdots, s_n), \underbar{t}=(t_1,t_2,
\cdots, t_n)\}$ is total in $\Gamma.$ \blema Let  $
U_{\underbar{s},\underbar{t}}(\underbar{u},\underbar{v}) \Omega, ~~
U_{\underbar{s}^\prime,\underbar{t}^\prime}(\underbar{p},\underbar{w})
\Omega \in \mathcal S.$\\
  Then there exist an integer $m\ge 1,~\underbar{a}=(a_1,a_2,
\cdots, a_m), \underbar{b}=(b_1,b_2, \cdots, b_n):~ 0 \le a_1\le
b_1\le  a_2  \le \ldots \le a_n\le  b_m< \infty,$ an  ordered
partition $R_1\cup R_2 \cup R_3=\{1,2,\cdots, m\}$ with $|R_i|=m_i$
and a  family of product vectors
$\underbar{x}_{k_l},\underbar{y}_{k_l} \in \mathbf h^{\otimes
m_1+m_2}, k_l\ge 1 :l=1,2\cdots,
m_1+m_2,~~\underbar{g}_{k_l},\underbar{h}_{k_l} \in \mathbf
h^{\otimes m_2+m_3}, k_l\ge 1 :l=1,2,\cdots, m_2+m_3$ such that \be
 U_{\underbar{s},\underbar{t}}(\underbar{u},\underbar{v})
=\sum_{\underbar{k}} \prod_{l\in R_1 \cup R_2} U_{a_l,b_l}
(x_{k_l},y_{k_l}) \ee

\be
 U_{\underbar{s}^\prime,\underbar{t}^\prime}(\underbar{p},\underbar{w})
=\sum_{\underbar{k}} \prod_{l\in R_2 \cup R_3} U_{a_l,b_l}
(g_{k_l},h_{k_l}). \ee \elema
\begin{proof}
 This follows  from the evolution hypothesis of  the  family  of
unitary  operators  $\{U_{s,t}\}$.
\end{proof}
\brmrk \label{V} Since the family of unitaries  $\{  V_{s,t}\}$ on
$\mathbf h \otimes \Gamma,$ enjoy all the properties  satisfied  by
the  family of unitaries  $\{  U_{s,t}\}$ on $\mathbf h \otimes
\mathcal H,$ the above Lemma also hold  if we replace $U_{s,t}$ by
$V_{s,t}.$

\ermrk

\blema For   $
U_{\underbar{s},\underbar{t}}(\underbar{u},\underbar{v}) \Omega, ~~
U_{\underbar{s}^\prime,\underbar{t}^\prime}(\underbar{p},\underbar{w})
\Omega \in \mathcal S.$\\
\be \label{uvinn} \langle
U_{\underbar{s},\underbar{t}}(\underbar{u},\underbar{v}) \Omega,
U_{\underbar{s}^\prime,\underbar{t}^\prime}(\underbar{p},\underbar{w})\Omega
\rangle=\langle
V_{\underbar{s},\underbar{t}}(\underbar{u},\underbar{v})
\textbf{e}(0),
V_{\underbar{s}^\prime,\underbar{t}^\prime}(\underbar{p},\underbar{w})
\textbf{e}(0) \rangle.
 \ee
  \elema
\begin{proof}
We have by previous Lemma  and {\bf Assumption: A}
 \bean &&\langle
U_{\underbar{s},\underbar{t}}(\underbar{u},\underbar{v}) \Omega,
U_{\underbar{s}^\prime,\underbar{t}^\prime}(\underbar{p},\underbar{w})\Omega
\rangle\\
&&=\sum_{\underbar{k}} \prod_{l\in R_1} \langle U_{b_l-a_l}
(x_{k_l},y_{k_l})\Omega ,\Omega\rangle \prod_{l\in R_2} \langle
U_{b_l-a_l} (x_{k_l},y_{k_l})\Omega
 U_{b_l-a_l} (g_{k_l},h_{k_l}) \Omega \rangle \\
 &&~~~~~~~~~~
 \prod_{l\in R_3} \langle
\Omega , U_{b_l-a_l} (g_{k_l},h_{k_l}) \Omega \rangle \\
&&=\sum_{\underbar{k}} \prod_{l\in R_1} \langle T_{b_l-a_l}
y_{k_l},x_{k_l}\rangle \prod_{l\in R_2} \langle g_{k_l} ,
 Z_{b_l-a_l} (| h_{k_l} >< y_{k_l} |)~x_{k_l} \rangle
 \prod_{l\in R_3} \langle
g_{k_l}  , T_{b_l-a_l} h_{k_l}\rangle \\
&&=\sum_{\underbar{k}} \prod_{l\in R_1} \langle V_{b_l-a_l}
(x_{k_l},y_{k_l})\textbf{e}(0), ~\textbf{e}(0)\rangle \prod_{l\in
R_2} \langle V_{b_l-a_l} (x_{k_l},y_{k_l})\textbf{e}(0)
 V_{b_l-a_l} (g_{k_l},h_{k_l})\textbf{e}(0) \rangle\\
 &&~~~~~~~~~~
 \prod_{l\in R_3} \langle
\textbf{e}(0) , V_{b_l-a_l} (g_{k_l},h_{k_l}) \textbf{e}(0)\rangle.
\eean Now  by Remark (\ref{V}), the above quantity is equal to
  $\langle
V_{\underbar{s},\underbar{t}}(\underbar{u},\underbar{v})
\textbf{e}(0),
V_{\underbar{s}^\prime,\underbar{t}^\prime}(\underbar{p},\underbar{w})
\textbf{e}(0) \rangle.$

\end{proof}

\bthm  There exist a unitary isomorphism $\Xi:\mathbf h \otimes
\mathcal H \rightarrow \mathbf h \otimes \Gamma$  such that \be
\label{U=V} U_t =\Xi^*~ V_t~ \Xi,~\forall~~ t\ge 0.\ee \ethm
\begin{proof}
Let us define a map $\Xi: \mathcal H \rightarrow \Gamma$ by setting,
for any $\xi=
U_{\underbar{s},\underbar{t}}(\underbar{u},\underbar{v}) \Omega \in
\mathcal S, ~~ \Xi
\xi:=V_{\underbar{s},\underbar{t}}(\underbar{u},\underbar{v})
\textbf{e}(0)\in \mathcal S^\prime$ and then extending linearly. So
by definition and totality of $\mathcal S^\prime,$ the range of
$\Xi$ is dense in $\Gamma.$
 To see that $\Xi$ is   a unitary  operator from $\mathcal H$ to
 $\Gamma$ it is enough to note  that
\be \langle   \Xi~ \xi, \Xi~ \xi^\prime  \rangle=\langle   \xi,
\xi^\prime  \rangle, ~ \forall~~ \xi,\xi^\prime \in \mathcal S\ee
which is  already proved in the  previous Lemma.

\noindent  Now  consider  the ampliated unitary operator $1_{\mathbf
h} \otimes \Xi $  from $\mathbf h \otimes \mathcal H $ to $\mathbf h
\otimes \Gamma$ and denote it by the same symbol $\Xi.$ In order to
prove (\ref{U=V})  it is enough to show that \be \langle   u\otimes
\xi , U_t v\otimes \xi^\prime  \rangle=\langle   \Xi(u\otimes  \xi),
V_t \Xi( v\otimes \xi^\prime ) \rangle, \forall~~ u,v\in \mathbf h,
\xi,\xi^\prime\in \mathcal S .\ee Note that
 $\Xi~U_t(u,v) \xi^\prime=V_t(u,v)~ \Xi~\xi^\prime.$  Now
by unitarity of $\Xi$, we have \bean &&  \langle   u\otimes  \xi ,
U_t v\otimes \xi^\prime  \rangle =\langle  \xi ,  U_t(u, v)
\xi^\prime \rangle
=\langle  \Xi ~\xi , \Xi~ U_t(u, v)   \xi^\prime  \rangle\\
&& =\langle  \Xi ~\xi ,  V_t(u, v)~  \Xi~ \xi^\prime \rangle=\langle
u\otimes  \Xi ~\xi,  V_t~  v\otimes  \Xi~ \xi^\prime \rangle=\langle
\Xi(u\otimes  \xi),  V_t~ \Xi( v\otimes \xi^\prime ) \rangle \eean
\end{proof}

\noindent {\bf Acknowledgements:} First author acknowledges National
Board for Higher Mathematics, DAE, India for  research support. Last
author acknowledges the partial support of Bhatnagar Fellowship
Project of CSIR, India. A portion of this work is completed at Delhi
Centre  of Indian Statistical  Institute.  All three authors have
benefited from the DST-DAAD (Indo-German) Exchange Programme at the
early stage of this work.

\end{document}